\documentclass{amsart}
\usepackage{amssymb}

\newtheorem{theorem}{Theorem}[section]
\newtheorem{proposition}[theorem]{Proposition}
\newtheorem{lemma}[theorem]{Lemma}
\newtheorem{corollary}[theorem]{Corollary}

\theoremstyle{definition}
\newtheorem*{definition}{Definition}

\theoremstyle{remark}
\newtheorem{remark}{Remark}[section]
\newtheorem*{note}{Note}

\numberwithin{equation}{section}

\begin{document}
\rightline{\scriptsize CRM-2369}

\title{Multivariable $q$-Racah polynomials}

\author{J. F. van Diejen}
\address{Centre de Recherches Math\'ematiques,
Universit\'e de Montr\'eal, C.P. 6128, succursale Centre-ville,
Montr\'eal (Qu\'ebec), H3C 3J7 Canada}

\thanks{The work of JFvD was supported in part by NSERC of Canada.}

\author{J. V. Stokman}
\address{Department of Mathematics, University of Amsterdam,
Plantage Muidergracht 24, 1018 TV  Amsterdam, The Netherlands}

%    General info

\date{June, 1996}


\begin{abstract}
The Koornwinder-Macdonald multivariable generalization
of the Askey-Wilson polynomials is studied for
parameters satisfying a truncation condition
such that the orthogonality measure becomes discrete
with support on a finite grid.
For this parameter regime the polynomials may be seen
as a multivariable counterpart of the (one-variable) $q$-Racah polynomials.
We present the discrete orthogonality measure, expressions
for the normalization constants converting the polynomials into
an orthonormal system (in terms of
the normalization constant for the unit polynomial),
and we discuss the limit $q\rightarrow 1$
leading to multivariable Racah type polynomials.
Of special interest is the situation that $q$ lies on the unit circle;
in that case it is found that there exists a natural parameter domain
for which the discrete orthogonality measure (which
is complex in general) becomes real-valued and positive.
We investigate the properties of a finite-dimensional
discrete integral transform for functions over the grid,
whose kernel is determined by the
multivariable $q$-Racah polynomials with parameters in this
positivity domain.
\end{abstract}

\maketitle

\section{Introduction}
Some years ago, Koornwinder \cite{koo:askey}
extended a construction of Macdonald \cite{mac:orthogonal}
(see also \cite{vil-klim:representations})
to arrive at a multivariable generalization of a family of
basic hypergeometric polynomials commonly known as the Askey-Wilson
polynomials \cite{ask-wil:some,gas-rah:basic}.
The multivariable polynomials of interest depend rationally on
a number of parameters and for parameter values in a certain
domain they form an orthogonal system with respect
to an explicitly given (positive)
continuous weight function with support on a
(real) $n$-dimensional torus (where $n$ denotes the number of variables).
Recently, it was shown that the parameter domain for
which the multivariable Askey-Wilson polynomials admit such
an interpretation
as orthogonal polynomials may be extended
if one allows the corresponding
orthogonality measure to have a partly continuous
and partly discrete support
\cite{sto:multivariable}. (Thus further generalizing the
corresponding situation in the case of one single variable, where the
phenomenon of discrete masses emerging in the Askey-Wilson
orthogonality measure was already known to occur \cite{ask-wil:some}.)

In the present paper we will demonstrate that for a different
parameter regime satisfying a certain truncation condition,
the multivariable Askey-Wilson polynomials can be
reduced to a finite-dimensional orthogonal system
with respect to a purely discrete weight function living on (i.e.
supported
on) a finite grid. The polynomials amount for these parameters
to a multivariable generalization of the $q$-Racah
polynomials \cite{ask-wil:set,gas-rah:basic}.
We will compute the normalization constants turning the polynomials
into an orthonormal system with
respect to the discrete orthogonality measure (in terms
of the corresponding normalization constant for the unit polynomial),
and also discuss the limit $q\rightarrow 1$
giving rise to multivariable Racah type polynomials.

In general we will not worry much about
the positivity of the weight function
and we will, in fact, for most of the time
allow parameters to be such that the discrete orthogonality
measure becomes complex. However, as will be outlined below
in further detail,
it is possible to restrict the parameters for $|q|=1$ in a
rather natural way to a subdomain such that the discrete orthogonality
measure for the multivariable $q$-Racah
polynomials becomes a genuine positive measure.
For parameters in this positivity domain the restriction of
the multivariable $q$-Racah polynomials to the grid points
entails (after renormalizing) an orthonormal basis
for the finite-dimensional Hilbert
space consisting of all (complex) functions over
the grid. From a functional-analytic
point of view, the renormalized polynomials
determine the kernel of a unitary finite-dimensional
integral transformation for these grid functions.

The material is structured as follows.
In Section~\ref{sec2} we first recall the definition of the
multivariable Askey-Wilson polynomials for generic parameters
as eigenfunctions of the Koornwinder-Macdonald
second order analytic $q$-difference operator.
Next, the discrete orthogonality measure for the polynomials is
introduced in Section~\ref{sec3} and it is explained how
the multivariable Askey-Wilson polynomials descend to a
$q$-Racah type finite-dimensional
orthogonal system when being restricted to the finite grid.
Crucial in the orthogonality proof is
the observation that the Koornwinder-Macdonald
second order $q$-difference operator
diagonalized by the polynomials is
symmetric with respect to the discrete inner product (just as
this observation turned out to be
essential when dealing with a purely
continuous or mixed continuous/discrete
orthogonality measure \cite{mac:orthogonal,koo:askey,sto:multivariable}).
The proof for the symmetry property of the $q$-difference operator
(in the discrete context)
is relegated to an appendix at the end of the paper (Appendix~\ref{appa}).
The orthonormalization constants
are given in Section~\ref{sec4}. Their computation, which is outlined
in Appendix~\ref{appb},
makes use of a recently introduced
system of recurrence relations (or Pieri type formulas) for the multivariable
Askey-Wilson polynomials \cite{die:self-dual,die:properties}.
In Section~\ref{sec5} we continue by
discussing the transition from the {\em basic} hypergeometric level
to the hypergeometric level ($q\rightarrow 1$).
In this limit---which one might also interpret as a transition from
trigonometric polynomials to rational counterparts obtained
by sending the period of the trigonometric functions
to infinity---our multivariable $q$-Racah type polynomials
degenerate into multivariable Racah type polynomials.
The paper is concluded
in Section~\ref{sec6} with the characterization of
a parameter domain
for which the weights determining the orthogonality
measure for the multivariable $q$-Racah polynomials become positive
when $q$ lies on the unit circle, and
a description of some properties of the finite-dimensional
discrete integral transform for grid functions
that is associated to the polynomials with parameters in this positivity
domain.
\begin{note}
Most objects of interest in this paper
(such as the polynomials, weight functions, difference equations,
normalization formulas, etc.)
depend rationally on a number of parameters.
Below it will always be assumed (unless explicitly stated
otherwise)
that the parameters are such that one stays away from
singularities (this is the case generically),
without repeatedly stressing this point each time.

We will be mainly concerned with the study of complex
orthogonality properties of our polynomials.
We shall say that a basis $\{ p_i\}_i$ for a
complex vector space $\mathcal{H}$
endowed with a nondegenerate bilinear
form $\langle \cdot ,\cdot\rangle$
is orthogonal
(with respect to $\langle \cdot ,\cdot\rangle$) if
$\langle p_i ,p_j\rangle =0$
for $i\neq j$. Furthermore, an orthogonal
basis of $\mathcal{H}$ will be called orthonormal (with
respect to $\langle \cdot ,\cdot\rangle$) if in addition
the normalization is such that the quantities
$\langle p_i ,p_i\rangle $ are all equal to one.
\end{note}

\section{Multivariable Askey-Wilson polynomials}\label{sec2}
In this section the definition of the multivariable
Askey-Wilson polynomials as eigenfunctions of the second
order Koornwinder-Macdonald $q$-difference operator is recalled.
The approach followed here is very much in the spirit of
Macdonald's original treatment in \cite{mac:orthogonal}.

Let ${\mathcal H}=\mathbb{C}[z_1,\ldots ,z_n,z_1^{-1},\ldots,z_n^{-1}]$ be the
space of Laurent polynomials in the variables $z_1,\ldots ,z_n$.
On this space the ($BC$ type) Weyl
group $W=S_n\ltimes (\mathbb{Z}_2)^n$ (i.e. the semidirect
product of the permutation group $S_n$ and the $n$-fold product
of the cyclic group $\mathbb{Z}_2$) acts naturally by
permutation and
inversion ($z_j\rightarrow z_j^{-1}$) of the variables
$z_1,\ldots ,z_n$.
The subspace ${\mathcal H}^W$ of ${\mathcal H}$
consisting of the $W$-invariant (Laurent) polynomials is
spanned by the symmetrized monomials
\begin{equation}\label{monom}
m_\lambda (z) = \sum_{\mu\in W (\lambda )} z_1^{\mu_1}\cdots z_n^{\mu_n},
\;\;\;\;\;\;\;\;\;\; \lambda\in\Lambda ,
\end{equation}
where $\Lambda(\cong  W\setminus \mathbb{Z}^n)$ denotes the integral cone
\begin{equation}\label{cone}
\Lambda = \{ \lambda\in \mathbb{Z}^n\; |\; \lambda_1\geq\lambda_2\geq
\cdots \geq \lambda_n\geq 0\}
\end{equation}
and the summation in \eqref{monom} is meant over the orbit of $\lambda$
under the action of the group $W$ (which
acts on vectors $\lambda\in\mathbb{Z}^n$ by
permuting and flipping the signs of the vector components
$\lambda_1,\ldots ,\lambda_n$). For future reference
we will partially order the
basis elements $m_\lambda$, $\lambda\in\Lambda$, by defining for
$\mu,\lambda\in\Lambda$ \eqref{cone}
\begin{equation}\label{order}
\mu\leq\lambda \;\;\; \text{iff} \;\;\;
\sum_{1\leq j\leq k}\mu_j \leq \sum_{1\leq j\leq k}\lambda_j \;\;\;
\text{for} \;\;\; k= 1,\ldots,n
\end{equation}
(and $\mu <\lambda$ if $\mu\leq\lambda$ with $\mu\neq\lambda$).

In \cite{koo:askey},
Koornwinder introduced the following generalization
of the  second order ($BC$ type) Macdonald $q$-difference operator
\begin{equation}\label{DAW}
D = \sum_{1\leq j\leq n}\Bigl( V_{j}(z)\, (T_{j,q}-1) \;+\;
V_{-j}(z)\, (T_{j,q}^{-1} -1)\Bigr) ,
\end{equation}
with
\begin{eqnarray*}
V_{\varepsilon j}(z)&=&
\frac{(1-t_0 z_j^{\varepsilon})
      (1-t_1 z_j^{\varepsilon})
      (1-t_2 z_j^{\varepsilon})
      (1-t_3 z_j^{\varepsilon})}
     {(1-z_j^{2\varepsilon})(1-q z_j^{2\varepsilon})} \\
& &\makebox[2em]{}\times \prod_{\stackrel{1\leq k\leq n}{k\neq j}}
\frac{(1-t\, z_j^{\varepsilon}\, z_k)(1-t\, z_j^{\varepsilon}\, z_k^{-1})}
     {(1-z_j^{\varepsilon}\, z_k)(1-z_j^{\varepsilon}\,z_k^{-1}) },
   \;\;\;\;\;\;\; \varepsilon =\pm 1, \\
\left(T_{j,q}f\right)(z)&=& f(z_1,\ldots,z_{j-1},q z_j,z_{j+1},
\ldots,z_n),
\end{eqnarray*}
and showed that this operator is triangular with respect to the
partially ordered basis of monomial symmetric functions:
\begin{equation}\label{tria}
D m_{\lambda} = \sum_{\mu\in\Lambda,\;\mu\leq\lambda}
E_{\lambda,\mu}\, m_{\mu}\;\;\;\;\text{with}\;\;\;\;
E_{\lambda,\mu}\in \mathbb{C}[q^{\pm 1},t,t_0,t_1,t_2,t_3]
\end{equation}
(i.e. the expansion coefficients (or matrix elements)
$E_{\lambda,\mu}$ depend polynomially on
$q^{\pm 1}$, $t$
and $t_0,\ldots ,t_3$).
The leading coefficient (or diagonal matrix element)
$E_{\lambda,\lambda}$ in \eqref{tria}
reads explicitly \cite{koo:askey}
\begin{equation}\label{eigenvAW}
E_{\lambda,\lambda}= \sum_{1\leq j\leq n}
\left(q^{-1}t_0t_1t_2t_3\, t^{2n-j-1}(q^{\lambda_j}-1)
\; +\; t^{j-1}(q^{-\lambda_j}-1)\right) .
\end{equation}
The triangularity \eqref{tria} of $D$
implies that the eigenvalue problem for the $q$-difference operator in
the space $\mathcal{H}^W$ is essentially finite-dimensional, because
it can be reduced to the invariant subspaces of the form
$\mathcal{H}^W_\lambda =
\text{Span} \{ m_\mu \}_{\mu\in\Lambda ,\mu \leq\lambda}$
(with $\lambda\in\Lambda$ \eqref{cone}).
The eigenvalues of $D$ in $\mathcal{H}^W$ are given by
the diagonal matrix elements $E_{\lambda,\lambda}$ \eqref{eigenvAW},
$\lambda\in\Lambda$. It is immediate from the explicit expression
in \eqref{eigenvAW}
that $E_{\mu ,\mu}\neq E_{\lambda ,\lambda}$ as a polynomial
in the parameters $q^{\pm 1}$, $t$ and $t_0,\ldots ,t_3$ if $\mu \neq \lambda$.
In other words, the eigenvalues are nondegenerate as
(Laurent) polynomial functions of the parameters.
The Koornwinder-Macdonald multivariable
Askey-Wilson polynomials are now---by definition---the corresponding
eigenfunctions $p_\lambda$, $\lambda\in\Lambda$.
\begin{definition}
The {\em multivariable Askey-Wilson polynomial}
associated with a (dominant weight) vector $\lambda\in \Lambda$ \eqref{cone}
is the (unique) monic
$W$-invariant Laurent polynomial of the form
\begin{subequations}
\begin{equation}\label{defcon1}
p_\lambda(z)=m_{\lambda}(z)+
       \sum_{\mu\in\Lambda ,\;\mu<\lambda}c_{\lambda,\mu}m_{\mu}(z)
\;\; \text{with}\;\; c_{\lambda,\mu}\in {\mathbb C}(q,t,t_0,t_1,t_2 ,t_3),
\end{equation}
such that
\begin{equation}\label{defcon2}
D\, p_\lambda=E_{\lambda,\lambda}\, p_\lambda.
\end{equation}
\end{subequations}
\end{definition}

Following Macdonald \cite{mac:orthogonal} (see also
\cite{die:self-dual,die:properties,sto-koo:limit})
it is possible to write down a somewhat more constructive formula
for the polynomial
$p_\lambda$ in terms of the monomial $m_\lambda$,
the operator $D$ and the eigenvalues
$E_{\mu ,\mu}$ ($\mu\leq\lambda$)
\begin{equation}\label{Macfor}
p_\lambda = \Biggl( \;\prod_{\mu\in\Lambda,\; \mu<\lambda}
\frac{D-E_{\mu,\mu}}{E_{\lambda,\lambda}-E_{\mu,\mu}}\Biggr)\; m_{\lambda} .
\end{equation}
(Notice that the r.h.s. is
well-defined as a rational expression
in the parameters in view of the fact that the denominators
are nonzero as (Laurent) polynomials in the parameters.)
The validity of this representation for
the multivariable Askey-Wilson polynomials is easily verified by inferring
that $p_\lambda$ \eqref{Macfor} satisfies the defining properties
\eqref{defcon1} and \eqref{defcon2}.
That the r.h.s. of \eqref{Macfor}
is of the form in \eqref{defcon1} is clear from the triangularity
of $D$; that it also satisfies \eqref{defcon2} is immediate from the fact that
the operator $\prod_{\mu\in\Lambda,\; \mu\leq \lambda} (D-E_{\mu ,\mu})$
annihilates the subspace
$\mathcal{H}^W_\lambda =
\text{Span} \{ m_\mu \}_{\mu\in\Lambda ,\mu \leq \lambda}$
as consequence of the Cayley-Hamilton theorem (and hence applying
$(D-E_{\lambda ,\lambda})$ to the r.h.s. of \eqref{Macfor} yields zero,
which is precisely \eqref{defcon2}).

In the one variable case, there is no $t$-dependence
and the eigenvalue equation $D\, p_\lambda = E_{\lambda ,\lambda}\, p_\lambda$
amounts in that case to the second order $q$-difference equation for the
Askey-Wilson polynomials.
Thus, the polynomials $p_\lambda$ then reduce to (monic)
Askey-Wilson polynomials \cite{ask-wil:some,gas-rah:basic}
\begin{equation}\label{bhrepAW}
p_\lambda (z) =
\frac{(t_0 t_1, t_0 t_2, t_0 t_3; q)_\lambda }{t_0^\lambda
(t_0 t_1 t_2 t_3 q^{\lambda -1};q)_\lambda }\;
{}_4\phi_3
\left(
\begin{matrix}
q^{-\lambda },\; t_0 t_1 t_2 t_3 q^{\lambda -1},\;
t_0 z,\; t_0 z^{-1} \\ [0.5ex]
t_0 t_1 ,\; t_0 t_2 ,\; t_0 t_3
\end{matrix} \; ; q,q \right) ,
\end{equation}
with $\lambda = 0,1,2,\ldots$
Here we have employed
standard notation (see e.g. \cite{ask-wil:some,gas-rah:basic})
for the basic hypergeometric series
\begin{equation*}
{}_{s+1}\phi_s\left( \begin{array}{c}
             a_1,\ldots ,a_{s+1} \\ b_1,\ldots ,b_s
            \end{array} ; q,z \right) =
\sum_{m=0}^\infty
\frac{(a_1,\ldots ,a_{s+1};q)_m}{(b_1,\ldots ,b_s;q)_m}
\frac{z^m}{(q;q)_m}
\end{equation*}
and the $q$-shifted factorials
\begin{eqnarray*}
(a_1,\ldots ,a_s;q)_m &=& (a_1;q)_m\cdots (a_s;q)_m, \\
(a;q)_m &=& (1-a)(1-aq)\cdots (1-aq^{m-1})
\end{eqnarray*}
(with $(a;q)_0\equiv 1$).

\section{A Discrete orthogonality measure}\label{sec3}
It follows from \cite{koo:askey} that if the parameters satisfy
the constraints
\begin{equation}\label{orthodomainAW}
0< q < 1,\;\;\; -1< t \leq 1, \;\;\; |t_r|\leq 1,\;\; r=0,1,2,3,
\end{equation}
with possible non-real parameters $t_r$ occurring in complex conjugate pairs
and pairwise products of the $t_r$ not lying in the interval $[1,\infty)$,
then the polynomials $p_\lambda$, $\lambda\in\Lambda$ constitute an orthogonal
system with respect to a continuous weight function
$\Delta^{\text{\tiny AW}}$:
\begin{equation}\label{orthogonalAW}
\int_{-\pi}^{\pi}\!\!\!\!\!\!\cdots
               \int_{-\pi}^{\pi}
p_\lambda (e^{ix})\, \overline{p_{\mu}(e^{ix})}\,
\Delta^{\text{\tiny AW}}  (x)\,
dx_1\cdots dx_n = 0 \;\;\; \text{if} \;\;\; \lambda \neq \mu
\end{equation}
where
\begin{eqnarray}\label{weightAW}
\Delta^{\text{\tiny AW}} (x) &=&
\prod_{\stackrel{1\leq j< k \leq n}
                {\varepsilon_1 ,\varepsilon_2 =\pm 1}}
\frac{(e^{i(\varepsilon_1 x_j
           +\varepsilon_2 x_k)}; q)_{\infty}}
     {(t\, e^{i(\varepsilon_1 x_j
           +\varepsilon_2 x_{k})}; q)_{\infty}} \\
& &  \times \prod_{\stackrel{1\leq j\leq n}{\varepsilon =\pm 1}}
\frac{( e^{2i \varepsilon x_j}; q)_\infty    }
     {(t_0\, e^{i \varepsilon x_j},\,
       t_1\, e^{i \varepsilon x_j},\,
       t_2\, e^{i \varepsilon x_j},\,
       t_3\, e^{i \varepsilon x_j};\, q)_\infty    } .  \nonumber
\end{eqnarray}
Here we have used the notation
$e^{ix}\equiv (e^{ix_1},\ldots ,e^{ix_n})$ and (cf. above)
\begin{equation*}
(a_1,\ldots ,a_r;q)_\infty = (a_1;q)_\infty \cdots (a_r;q)_\infty,
\;\;\;\;
(a;q)_\infty = \prod_{m=0}^\infty (1-aq^m).
\end{equation*}
It should be noted that at this point of our
presentation it is not even obvious that the
multivariable Askey-Wilson polynomials
defined in the previous section actually exist for {\em all} parameter
values in the above domain, since the domain might a priori contain
parameter values for which
the coefficients $c_{\lambda ,\mu}$ \eqref{defcon1} have
a singularity.
However, that such singularities indeed do not occur
is seen from an alternative characterization of the
multivariable Askey-Wilson polynomial $p_\lambda$
as the polynomial of the form
\begin{subequations}
\begin{equation}\label{altcon1}
p_\lambda(z)=m_{\lambda}(z)+
       \sum_{\mu\in\Lambda ,\;\mu<\lambda}c_{\lambda,\mu}m_{\mu}(z)
\;\; \text{with}\;\; c_{\lambda,\mu}\in {\mathbb C},
\end{equation}
satisfying
\begin{equation}\label{altcon2}
\int_{-\pi}^{\pi}\!\!\!\!\!\!\cdots
               \int_{-\pi}^{\pi}
p_\lambda (e^{ix})\, \overline{m_{\mu}(e^{ix})}\,  \Delta^{\text{\tiny AW}}
(x)\,
dx_1\cdots dx_n = 0 \;\;\; \text{for} \;\;\; \mu < \lambda .
\end{equation}
\end{subequations}
(In \cite{koo:askey} these two properties were in fact used
to {\em define} the multivariable Askey-Wilson polynomials.)
It is clear that the polynomials determined
by the conditions \eqref{altcon1}, \eqref{altcon2}
are well-defined for all parameter values in the domain given
by \eqref{orthodomainAW} and, furthermore, that they are
continuous in $t$ and $t_0,\ldots ,t_3$ for these parameter values.
By showing (as was done in \cite{koo:askey})
that the difference operator $D$ \eqref{DAW}
is symmetric with respect
to the $L^2$ inner product with weight function $\Delta^{\text{\tiny AW}}$,
i.e., that
\begin{eqnarray*}
&&\int_{-\pi}^{\pi}\cdots \int_{-\pi}^{\pi}
 (D m_\lambda) (e^{ix})\,  \overline{m_\mu (e^{ix})}
\Delta^{\text{\tiny AW}}(x)
dx_1\dots dx_n =\\
&&\makebox[3em]{}
\int_{-\pi}^{\pi}\cdots\int_{-\pi}^{\pi} m_\lambda (e^{ix})\, \overline{(Dm_\mu
)(e^{ix})}
\Delta^{\text{\tiny AW}}(x)dx_1\cdots dx_n
\end{eqnarray*}
and combining this with the triangularity property in
\eqref{tria}, one finds that the polynomial $p_\lambda (z)$  of
the form \eqref{altcon1}, \eqref{altcon2} satisfies
the eigenvalue equation
$D\, p_\lambda = E_{\lambda ,\lambda}\, p_\lambda$, which shows that
for the parameter domain of interest this alternative characterization
leads us to the same polynomial $p_\lambda$ as in the previous section
(and hence that the coefficients $c_{\lambda ,\mu}$ are regular
for these parameter values).
The orthogonality \eqref{orthogonalAW} of the multivariable Askey-Wilson
polynomials $p_\lambda$, $\lambda\in\Lambda$
with respect to
the weight function $\Delta^{\text{\tiny AW}}$ \eqref{weightAW}
now follows, first
for generic parameters in the above domain from the fact that they are
eigenfunctions of a symmetric operator $D$ corresponding to different
(real) eigenvalues, and then for all parameters in this domain
by a continuity argument \cite{koo:askey}.

Very recently, it was shown
that the parameter domain for which the multivariable Askey-Wilson
polynomials admit an interpretation as orthogonal polynomials may
be further extended.
Specifically, it follows from
\cite{sto:multivariable} that for $t=q^m$, with $m$ an arbitrary
nonnegative integer,
one can remove the constraints $|t_r|\leq 1$ in \eqref{orthodomainAW}
(while still keeping the other
restrictions on these parameters though) to end up with an orthogonality
relation for the polynomials $p_\lambda$, $\lambda\in\Lambda$ consisting of
a continuous part of the form \eqref{orthogonalAW} and additional mixed
continuous,
mixed discrete parts.

The main purpose of the present section is to demonstrate
that for parameters satisfying a truncation condition
the multivariable Askey-Wilson polynomials
give rise to a finite-dimensional orthogonal system
with a purely discrete weight function living
on a finite grid.

First it is needed to introduce some notation.
Let $N$ be an arbitrary nonnegative integer and let
$\Lambda_{N}$
be the alcove of dominant weight vectors of the form
\begin{equation}\label{alcove}
\Lambda_N = \{ \lambda\in \Lambda\; |\; \lambda \leq N\omega \}
\;\;\; \text{with}\;\;\; \omega \equiv e_1+\cdots +e_n.
\end{equation}
(Here and below $e_j$ represents the
$j$th unit vector in the standard basis of $\mathbb{R}^n$.)
We will show that for parameters
subject to a truncation condition of the type
\begin{equation}\label{trunc}
t_at_bt^{n-1}=q^{-N},
\end{equation}
with $a,b \in \{ 0,1,2,3\}$ (fixed) and $b\neq a$, the multivariable
Askey-Wilson polynomials $p_\lambda$, $\lambda\in\Lambda_N$
form a finite-dimensional orthogonal system
with respect to a (generally complex) discrete measure
with support on the grid points
\begin{equation}\label{grid}
\tau q^{\nu} \equiv (\tau_1q^{\nu_1},\ldots ,
         \tau_jq^{\nu_j},\ldots , \tau_n q^{\nu_n}) ,\;\;\;\;\;
\nu \in \Lambda_N
\end{equation}
where
\begin{equation}\label{tau}
\tau_j = t^{n-j} t_a, \;\;\;\;\;\;\;\;\;\;\;\; j=1,\ldots ,n.
\end{equation}
Specifically, the (complex) orthogonality relation becomes for these parameters
(i.e. satisfying condition \eqref{trunc})
\begin{equation}\label{orthogonalqR}
\sum_{\nu \in\Lambda_N }
p_\lambda (\tau q^\nu) p_{\mu} (\tau q^\nu) \Delta^{\text{\tiny qR}} (\nu) = 0
\;\;\; \text{for} \;\;\;\lambda \neq \mu
\;\;\;\;\;\;\; (\lambda ,\mu \in \Lambda_N),
\end{equation}
where the weights $\Delta^{\text{\tiny qR}} (\nu )$, $\nu\in\Lambda_N$
are given explicitly by
\begin{equation}\label{weightqR}
\Delta^{\text{\tiny qR}} (\nu ) =
\frac{ 1 }
{C_+^{\text{\tiny qR}}(\nu )\: C_-^{\text{\tiny qR}}(\nu ) } ,
\end{equation}
with
\begin{eqnarray}\label{Cqr+}
C_+^{\text{\tiny qR}}(\nu )\!\!\! &=&\!\!\! c_0(\nu)
\prod_{1\leq j< k\leq n}
\left( \frac{(\tau_j\tau_k;q)_{\nu_j+\nu_k}}{(t\tau_j\tau_k;q)_{\nu_j+\nu_k}}
\frac{(\tau_j\tau_k^{-1};q)_{\nu_j-\nu_k}}
     {(t\tau_j\tau_k^{-1};q)_{\nu_j-\nu_k}} \right) \\
&& \makebox[2em]{}\times \prod_{1\leq j\leq n}
\left( \frac{(\tau_j^2;q)_{2\nu_j}}
     {\prod_{0\leq r\leq 3}(t_r\tau_j;q)_{\nu_j}} \right) ,\nonumber\\
\label{Cqr-}
C_-^{\text{\tiny qR}}(\nu )\!\!\! &=&\!\!\! c_0(\nu)
\prod_{1\leq j< k\leq n}
\left( \frac{(t^{-1}q\tau_j\tau_k ;q)_{\nu_j+\nu_k}}
     {(q\tau_j\tau_k ;q)_{\nu_j+\nu_k}}
\frac{(t^{-1}q\tau_j\tau_k^{-1} ;q)_{\nu_j-\nu_k}}
     {(q\tau_j\tau_k^{-1} ;q)_{\nu_j-\nu_k}} \right) \\
& &\makebox[2em]{} \times \prod_{1\leq j\leq n}
\left( \frac{\prod_{0\leq r\leq 3}(t_r^{-1}q\tau_j;q)_{\nu_j}}
     {(q\tau_j^2;q)_{2\nu_j}} \right) \nonumber
\end{eqnarray}
and
\begin{equation}\label{c0}
c_0(\nu) = \prod_{1\leq j\leq n}
\left (t^{n-j}(t_0t_1t_2t_3q^{-1})^{1/2} \right)^{\nu_j} .
\end{equation}
Here we have introduced discretized
functions $C_\pm^{\text{\tiny qR}}(\nu )$ that
are reminiscent (upon dualization) of the $c$-functions of Harish-Chandra
in harmonic analysis (see e.g. \cite{hel:groups,hec-sch:harmonic}).
\begin{lemma}\label{nonzero}
a. For $\nu\in\Lambda$ \eqref{cone} the functions
$C_\pm^{\text{\tiny qR}}(\nu )$ \eqref{Cqr+}, \eqref{Cqr-}
are well-defined and nonzero as
meromorphic expressions in the parameters $t$, $t_0,\ldots ,t_3$
and $q$.

b. For $\nu\in\Lambda_N$ \eqref{alcove} the functions
$C_\pm^{\text{\tiny qR}}(\nu )$ \eqref{Cqr+}, \eqref{Cqr-}
are well-defined and nonzero as
meromorphic expressions in the parameters $t$, $t_0,\ldots ,t_3$
and $q$ subject to the truncation condition \eqref{trunc}.
\end{lemma}
Here and below when stating that an expression is meromorphic or rational
in the parameters $t$, $t_0,\ldots ,t_3$
and $q$ subject to the truncation condition \eqref{trunc}, it is meant
that after elimination of $t_b$ (or $t_a$) with the aid of
the relation $t_at_b = t^{1-n}q^{-N}$, the resulting
expression is meromorphic/rational in the remaining parameters
$t_a$ (or $t_b$), $t_c, t_d, t$ and $q$
(where $t_c$ and $t_d$ denote the two parameters
complementing $t_a$ and $t_b$ in $\{ t_0,t_1,t_2,t_3\}$
such that $\{ t_a,t_b,t_c,t_d\} =\{ t_0,t_1,t_2,t_3\}$).

The proof of this lemma is immediate from inspection of the explicit
expressions for $C_\pm^{\text{\tiny qR}}(\nu )$ given above
(no factor in the numerator or the
denominator of $C_\pm^{\text{\tiny qR}}(\nu )$ becomes identical to zero).
Observe also that the functions $C_\pm^{\text{\tiny qR}}(\nu )$
are rational in the parameters except for
the square roots appearing in the common
factor $c_0(\nu)$ \eqref{c0}.
In the discretized
weight function $\Delta^{\text{\tiny qR}} (\nu )$ \eqref{weightqR}
these square roots can be collected
(there is an even number of them) and rationality becomes restored.
Notice also that for parameters subject to the truncation condition
$C_+^{\text{\tiny qR}}(\nu )$ \eqref{Cqr+} would be infinite
if $\nu\in\Lambda\setminus\Lambda_N$ (i.e., when $\nu$ lies in
the cone $\Lambda$ \eqref{cone} but outside
the alcove $\Lambda_N$ \eqref{alcove}). In that case we have that
$\nu_1 > N$ and hence
the factor $(t_b\tau_1;q)_{\nu_1}=(t_at_bt^{n-1};q)_{\nu_1}$
in the denominator of $C_+^{\text{\tiny qR}}(\nu )$ \eqref{Cqr+}
is identical to zero when the truncation condition \eqref{trunc} holds.
As a result, the weight function $\Delta^{\text{\tiny qR}} (\nu )$
\eqref{weightqR} vanishes on $\Lambda\setminus\Lambda_N$ for
parameters subject to the truncation condition.

The orthogonality relation \eqref{orthogonalqR} should be read
as a set of identities for the multivariable Askey-Wilson polynomials
$p_\lambda (z)$, $\lambda \in \Lambda_N$ with parameters subject
to the truncation condition \eqref{trunc}.
It is clear from Macdonald's representation in \eqref{Macfor} and the explicit
formula for the eigenvalues $E_{\lambda ,\lambda}$ \eqref{eigenvAW}
that the polynomials $p_\lambda (z)$ with this restriction on the
parameters are well-defined as
rational expressions in the parameters
subject to the truncation condition
(no denominator in \eqref{Macfor} becomes identical to zero
after imposing the truncation condition \eqref{trunc}).
Since the same is true for the weights
$\Delta^{\text{\tiny qR}}(\nu)$ \eqref{weightqR},
the terms in the l.h.s. of \eqref{orthogonalqR} make sense as
rational expressions in the parameters
$q$, $t$ and $t_0,t_1, t_2,t_3$ subject to condition~\eqref{trunc}.

In the case of one single variable ($n=1$) the orthogonality relation
\eqref{orthogonalqR} reduces to the well-known orthogonality relation
\begin{eqnarray}\label{orthoqR1}
\sum_{0\leq \nu \leq N}
p_\lambda (t_a q^{\nu}) p_{\mu} (t_a q^{\nu}) \Delta^{\text{\tiny qR}} (\nu )=
0 &&\\
\text{for}& & \!\!\!\!\!\!\!\lambda \neq \mu\;\;\;\;\;\;\;
(\lambda ,\mu \in \{ 0,\ldots ,N\} ),\nonumber
\end{eqnarray}
with
\begin{equation}\label{deltaqR}
\Delta^{\text{\tiny qR}} (\nu ) =
\frac{(1-t_a^2 q^{2\nu})}
     {(t_0 t_1 t_2 t_3 q^{-1})^{\nu} (1-t_a^2)}
\frac{(t_0t_a, t_1t_a, t_2t_a, t_3t_a;q)_{\nu}}
     {(t_0^{-1}qt_a, t_1^{-1}qt_a,
                         t_2^{-1}qt_a, t_3^{-1}qt_a;q)_{\nu}},
\end{equation}
for the (monic) Askey-Wilson polynomials given
by $p_\lambda (z)$ \eqref{bhrepAW}
with parameters subject to a condition of the type
$t_at_b=q^{-N}$ ($b\neq a$) \cite{ask-wil:set,gas-rah:basic}.
In the discrete case, with parameters subject to the truncation condition,
these one-variable polynomials are usually referred to in the literature
as $q$-Racah polynomials rather than Askey-Wilson polynomials.

In the above formulas (and throughout the paper) $t_a$ and $t_b$ denote
two fixed, distinct, but otherwise
arbitrary parameters from the set $\{ t_0,t_1,t_2,t_3\}$.
Since the monic (multivariable) Askey-Wilson polynomials
$p_\lambda (z)$ are symmetric with respect to permutations of
the parameters $t_0,\ldots ,t_3$ (this is immediate from their definition
as monic eigenfunctions of the difference operator $D$ \eqref{DAW}), the
actual choices of $t_a$ and $t_b$ are all equivalent up to permutation
of the parameters $t_0,t_1,t_2,t_3$.
Notice, however, that the choice of $t_a$ {\em does} influence
the position of the grid points $\tau q^\nu$ \eqref{grid}.
In the case of one single variable one usually works with
a different normalization of the $q$-Racah polynomials that
breaks the permutation symmetry and designates a preferred role
to the parameter $t_0$. (Instead of employing the monic polynomials
of \eqref{bhrepAW} one then
omits the constant factor in front of the terminating ${}_4\phi_3$ series.)
It is in that case custom to take $t_a$ to be $t_0$ (thus fixing
the grid points to $t_0 q^\nu$, $\nu =0,\ldots ,N$).
This then leaves three possible choices for the
other parameter $t_b$ entering the truncation condition $t_at_b=q^{-N}$,
which leads one to distinguish the three cases
$t_0t_1=q^{-N}$, $t_0t_2=q^{-N}$ and $t_0t_3=q^{-N}$.
Here instead of treating each of these cases separately we prefer not
to make such distinctions and have tried rather to
emphasize the symmetry in the parameters $t_0,\ldots ,t_3$ with our notation.

In order to prove the discrete orthogonality relations \eqref{orthogonalqR}
and analyze their properties in more detail we will need some further notation.
Specifically, let us introduce a bilinear form
$\langle \cdot ,\cdot \rangle_N^{\text{\tiny qR}}$
on the finite-dimensional subspace (of $\mathcal{H}^W$)
\begin{equation}\label{HN}
\mathcal{H}^{\text{\tiny qR}}_N \equiv \text{Span} \{ m_\lambda
\}_{\lambda\in\Lambda_N}
\end{equation}
determined by
\begin{equation}\label{discreteIP}
\langle f ,g \rangle_N^{\text{\tiny qR}} =
\sum_{\nu \in\Lambda_N} f (\tau q^{\nu})
 g (\tau q^{\nu}) \Delta^{\text{\tiny qR}} (\nu)
\;\;\;\;\;\;\;\;\;\;\; (f,g \in\mathcal{H}^{\text{\tiny qR}}_N) .
\end{equation}
The orthogonality of the multivariable Askey-Wilson (or $q$-Racah)
polynomials with respect to the discrete weight function
$\Delta^{\text{\tiny qR}}$ \eqref{weightqR}
can now be phrased in terms of the following theorem.
\begin{theorem}\label{orthothm}
For parameters $q$, $t$ and $t_0,\ldots ,t_3$ subject to the
truncation condition
\begin{equation*}
t_at_b t^{n-1}=q^{-N}\;\;\;\;\;\;\;\;\;\; b\neq a
\end{equation*}
(where $N$ denotes an arbitrary nonnegative integer and $a$, $b$ are
two fixed, distinct but otherwise arbitrary numbers
from the set $\{ 0,1,2,3\}$),
the multivariable Askey-Wilson polynomials
$p_\lambda (z)$, $\lambda\in\Lambda_N$ \eqref{alcove}
satisfy the orthogonality relation
\begin{equation}\label{orthoqR}
\langle p_\lambda ,p_{\mu} \rangle_N^{\text{\tiny qR}} = 0\;\;
\text{for}\;\; \lambda \neq \mu
\;\;\;\;\;\;\;\;\;\;\; (\lambda ,\mu\in\Lambda_N),
\end{equation}
with the bilinear form $\langle \cdot ,\cdot \rangle_N^{\text{\tiny qR}}$
being defined by \eqref{discreteIP}.
\end{theorem}
The proof of the orthogonality theorem operates along the same lines
as the orthogonality proof in \cite{koo:askey} for the continuous
case (cf. the beginning of this section) and is based on the
statement that (for parameters satisfying
the truncation condition)
the difference operator $D$ \eqref{DAW} is symmetric with
respect to the bilinear form
$\langle \cdot ,\cdot \rangle_N^{\text{\tiny qR}}$.
\begin{proposition}\label{symprp}
For parameters in accordance with the assumptions in
Theorem~\ref{orthothm},
the difference operator $D$ \eqref{DAW}
is symmetric with respect to the bilinear form
$\langle\cdot ,\cdot\rangle_N^{\text{\tiny qR}}$ \eqref{discreteIP}, i.e.
\begin{equation}
\langle D f , g\rangle_N^{\text{\tiny qR}} = \langle f
, D g \rangle_N^{\text{\tiny qR}}
\;\;\;\;\;\;\;\;\;\;\;\;\; (f , g\in\mathcal{H}^{\text{\tiny qR}}_N ).
\end{equation}
\end{proposition}
The proof of this proposition has, in turn, been relegated
to Appendix~\ref{appa} at the end of the paper.
The orthogonality \eqref{orthoqR} now follows from
the fact that
the polynomials $p_\lambda$, $\lambda\in\Lambda_N$
are eigenfunctions of a symmetric operator
$D$ corresponding
to eigenvalues $E_{\lambda ,\lambda}$ \eqref{eigenvAW}, $\lambda\in\Lambda_N$
that are nondegenerate (i.e. distinct for different $\lambda$)
as rational functions in the parameters $q$, $t$ and
$t_0,\ldots ,t_3$ subject to condition \eqref{trunc}.

An important property of the bilinear form
$\langle \cdot ,\cdot \rangle_N^{\text{\tiny qR}}$ is that for generic
parameters
(subject to the truncation condition)
it is nondegenerate on the space $\mathcal{H}^{\text{\tiny qR}}_N$ \eqref{HN}
(i.e., if for a certain $f\in\mathcal{H}^{\text{\tiny qR}}_N$
one has that $\langle f,g\rangle_N^{\text{\tiny qR}} =0$ for all $g\in
\mathcal{H}^{\text{\tiny qR}}_N$,
then $f$ must be zero).
\begin{proposition}\label{nondegprp}
The bilinear form $\langle \cdot ,\cdot \rangle_N^{\text{\tiny qR}}$
\eqref{discreteIP}
is nondegenerate
on the space $\mathcal{H}^{\text{\tiny qR}}_N$ \eqref{HN}
for generic parameters $q$, $t$ and
$t_0,\ldots ,t_3$ (subject to condition \eqref{trunc}).
\end{proposition}
The proof of Proposition~\ref{nondegprp} readily
follows after recalling that the weights $\Delta^{\text{\tiny qR}}(\nu)$,
$\nu\in \Lambda_N$ are nonzero as rational expressions
in the parameters subject to the truncation condition
(no weight becomes identical to zero) and
combining this with the following lemma.
\begin{lemma}\label{replem}
For generic  $q$, $t$ and $t_a$
every function in the space $\mathcal{H}^{\text{\tiny qR}}_N$
\eqref{HN}
is {\em uniquely} characterized by
its values on the grid points $\tau q^{\lambda}$,
$\lambda\in\Lambda_N$.
\end{lemma}
To prove Lemma~\ref{replem} one uses that
the determinant of the matrix
\begin{equation}\label{qmatrix}
M = [m_\mu (\tau q^{\nu})]_{\mu ,\nu\in\Lambda_N}
\end{equation}
(where
the columns and rows are ordered by means of some (any) total order
of the weight vectors in the cone $\Lambda$ \eqref{cone}
extending the partial order \eqref{order}) is nonzero as a
Laurent polynomial in $q$ for generic $t$ and $t_a$(this follows from
the dominant behavior of the determinant
for $q\rightarrow +\infty$, see \cite[Lemma~4]{die:self-dual}).
Hence, if $\sum_{\mu \in \Lambda_N} c_{\mu} m_{\mu}(\tau q^{\nu}) =0$
for all $\nu \in\Lambda_N$, then the coefficients $c_\mu$ must all be zero.

Since the dimension of the space of {\em all} functions over the grid
\eqref{grid} is equal to the number of points in the grid, which is for
generic $q$, $t$ and $t_a$
equal to the number of points in the alcove $\Lambda_N$
\eqref{alcove}
and hence identical to the dimension of the space
$\mathcal{H}^{\text{\tiny qR}}_N$ \eqref{HN},
Lemma~\ref{replem} actually says---when read in an appropriate
way---that for generic $q$, $t$ and $t_a$
{\em any} complex function defined on the grid points
$\tau q^{\lambda}$, $\lambda\in\Lambda_N$ may be
identified with the restriction of a function in the space
$\mathcal{H}^{\text{\tiny qR}}_N$ to the grid.
In other words, for such generic $q$, $t$ and $t_a$ every
function on the grid can be approximated {\em exactly}
by (the restriction to the grid of) a {\em unique} function in
$\mathcal{H}^{\text{\tiny qR}}_N$.
\begin{proposition}\label{repprp}
For generic $q$, $t$ and $t_a$
{\em any} complex function over the grid with points
$\tau q^\nu$, $\nu\in\Lambda_N$ can be represented
{\em exactly} by the restriction
to the grid points of a {\em unique} function in the space
$\mathcal{H}^{\text{\tiny qR}}_N$ \eqref{HN}.
\end{proposition}

In principle there is no real reason why the definition of
the bilinear form
$\langle \cdot ,\cdot \rangle_N^{\text{\tiny qR}}$ \eqref{discreteIP}
should be confined to the subspace $\mathcal{H}^{\text{\tiny qR}}_N$ and
could not, e.g., be extended
to the whole space $\mathcal{H}^W$ (which consist of all
$W$-invariant Laurent polynomials in $n$ variables).
In fact, this is why
the condition $f,g \in \mathcal{H}^{\text{\tiny qR}}_N$ in \eqref{discreteIP}
has been put between parenthesis.
It turns out that Theorem~\ref{orthothm} and
Proposition~\ref{symprp}
remain valid (with the given proofs applying verbatim)
if one allows $\lambda $ and $\mu$ to be any dominant weight
vector from the cone $\Lambda$ \eqref{cone} and
$f$ and $g$ to be any function
in the space $\mathcal{H}^W$, respectively.
However, since such a generalization of
Theorem~\ref{orthothm} would imply that
$p_\lambda $ is orthogonal to $\mathcal{H}^{\text{\tiny qR}}_N$
if $\lambda$ is not in $\Lambda_N$, it follows
from combination
of  Proposition~\ref{nondegprp} and Proposition~\ref{repprp}
that in that case $p_\lambda$ must
actually be zero on the grid.

\begin{proposition}\label{findimredprp}
For parameters subject to the truncation condition
\eqref{trunc} one has that
\begin{equation}
p_\lambda (\tau q^{\nu})=0\;\;\;\;\;
\text{for}\;\; \nu\in\Lambda_N
\end{equation}
if $\lambda\in\Lambda\setminus\Lambda_N$.
\end{proposition}
The upshot of Proposition~\ref{findimredprp} is that the multivariable
Askey-Wilson polynomials $p_\lambda (z)$, $\lambda\in\Lambda$ with
parameters subject to the truncation condition \eqref{trunc}
descend to a finite-dimensional orthogonal system
spanned by $p_\lambda (z)$, $\lambda\in\Lambda_N$
when being restricted
to the grid points $\tau q^{\nu}$, $\nu\in\Lambda_N$.
Furthermore, the results of this section can now be
alternatively summarized in the following
way: for generic $q$ and parameters subject to the truncation condition,
the restriction of the multivariable Askey-Wilson
polynomials $p_\lambda (z)$, $\lambda\in\Lambda_N$ to the grid points
$\tau q^{\nu}$, $\nu\in\Lambda_N$ yields an
orthogonal basis for the (finite-dimensional) space of functions
over this grid endowed with the bilinear form
$\sum_{\nu\in\Lambda_N} f(\tau q^{\nu}) g(\tau q^{\nu})
\Delta^{\text{\tiny qR}} (\nu)$.

\begin{remark}\label{cext}
It is convenient to think of the discrete weight function determined by
the weights $\Delta^{\text{\tiny qR}}(\nu )$
\eqref{weightqR} as arising from
a continuous weight function $\bar{\Delta}^{\text{\tiny qR}}(z)$ restricted
to the grid points $\tau q^\nu$ \eqref{grid} (such that
$\bar{\Delta}^{\text{\tiny qR}}(\tau q^\nu )=\Delta^{\text{\tiny qR}}(\nu )$).
In principle such a function
$\bar{\Delta}^{\text{\tiny qR}}(z)$
of course exists and is a priori highly non-unique.
For $0<|q|<1$ a possibility would be to take
\begin{eqnarray}
\bar{\Delta}^{\text{\tiny qR}}(z) &=&
\mathcal{D}_0
\prod_{1\leq j\leq n} (t^{2(n-j)}t_0t_1t_2t_3q^{-1})^{-\log (z_j)/\log (q)}
                          \nonumber  \\
& & \times \prod_{1\leq j< k\leq n} \left( (1-z_jz_k)(1-z_jz^{-1}_k)
\frac{(t^{-1}qz_jz_k,t^{-1}qz_jz_k^{-1};q)_\infty}
     {(t z_jz_k,t z_jz_k^{-1};q)_\infty} \right)  \label{context} \\
& & \times \prod_{1\leq j\leq n}
\left( (1-z_j^2)
\frac{(t_0^{-1}qz_j,t_1^{-1}qz_j,t_2^{-1}qz_j,t_3^{-1}qz_j;q)_\infty}
     {(t_0z_j,t_1 z_j,t_2 z_j,t_3 z_j;q)_\infty} \right) , \nonumber
\end{eqnarray}
where $\mathcal{D}_0$ denotes a normalization constant
determined by the requirement
that evaluation of $\bar{\Delta}^{\text{\tiny qR}}(z)$ in
$z=\tau\equiv (\tau_1,\dots ,\tau_j,\ldots ,\tau_n)$---i.e.,
$z=\tau q^\nu$ \eqref{grid} with $\nu =0$---should
yield the value one (since $\Delta^{\text{\tiny qR}}(0)=1$).
\end{remark}
\begin{remark}\label{gridstruct}
It is instructive to view the grid points \eqref{grid} on which
the discrete masses of the orthogonality measure
determined by $\Delta^{\text{\tiny qR}}$
are positioned as being of the form $q^{\rho+\nu}$, $\nu \in \Lambda_N$,
where $\rho$ denotes an $n$-dimensional
vector that is related to $\tau$ \eqref{tau} by a logarithm
\begin{equation}
\rho = (\rho_1,\ldots ,\rho_n)\;\;\;\text{with}\;\;
\rho_j =\log (\tau_j)/\log (q), \;\; j=1,\ldots, n.
\end{equation}
Thus, up to an exponential scaling ($z\rightarrow q^z$) the grid consists
of the points in the alcove $\Lambda_N$ \eqref{alcove} translated
over the vector $\rho$.
\end{remark}

\section{Orthonormalization}\label{sec4}
In this section we will present the
normalization constants
turning the multivariable Askey-Wilson
polynomials $p_\lambda (z)$,
$\lambda\in\Lambda_N$ \eqref{alcove}---with parameters subject
to the truncation condition \eqref{trunc}---into an orthonormal
basis for the space $\mathcal{H}^{\text{\tiny qR}}_N$ \eqref{HN}
with respect to the bilinear form
$\langle \cdot ,\cdot \rangle_N^{\text{\tiny qR}}$ \eqref{discreteIP}.
To this end it is needed to compute the quantities
$\langle p_\lambda ,p_\lambda\rangle_N^{\text{\tiny qR}}$,
$\lambda\in\Lambda_N$.
The evaluation of these sums pivots on a previously
introduced system of Pieri type recurrence formulas
for the multivariable Askey-Wilson polynomials
\cite{die:self-dual,die:properties}. At the time
of their introduction the
same formulas were also used to verify
the explicit expressions (conjectured by Macdonald)
for the squared norm of $p_\lambda(z)$
with respect to Koornwinder's inner product \eqref{orthogonalAW}.
In both cases (i.e. the discrete and the continuous case)
the mechanism leading to the solution of the orthonormalization problem
is very similar.
Hence, we will refrain from providing full details here and refer
the reader to Appendix~\ref{appb}, where
the main ingredients of the recipe are outlined
with an emphasis on those points at which the
discrete case differs from the continuous case.

We would like to add that exactly the same approach
may also serve to obtain the orthonormalization constants
(in terms of the norm of the unit polynomial)
for the parameter regime with the mixed continuous/discrete
orthogonality measure considered in \cite{sto:multivariable}.
The answer will in that case be formally identical
to that in the case of Koornwinder's purely continuous
orthogonality measure $\Delta^{\text{\tiny AW}}$ \eqref{weightAW},
except that the parameter domain now gets extended in the
way indicated at the beginning of the previous section.

In order to describe the evaluation formula for
$\langle p_\lambda , p_\lambda \rangle_N^{\text{\tiny qR}}$,
it is convenient to introduce dual parameters $\hat{t}_r$,
$r=0,\ldots ,3$ that are related to the parameters $t_r$,
$r=0,\ldots ,3$ in the following way
\begin{equation}\label{hattr}
\begin{array}{lll}
\hat{t}_a &=& (t_a t_b      t_c      t_d      q^{-1})^{1/2},\\
\hat{t}_b &=& (t_a t_b      t_c^{-1} t_d^{-1} q     )^{1/2},\\
\hat{t}_c &=& (t_a t_b^{-1} t_c      t_d^{-1} q     )^{1/2},\\
\hat{t}_d &=& (t_a t_b^{-1} t_c^{-1} t_d      q     )^{1/2},
\end{array}
\end{equation}
where $c$ and $d$ denote
the two indices that complement the
indices $a$ and $b$ entering the truncation
condition \eqref{trunc}
such that $\{ a,b,c,d\} =\{ 0,1,2,3\}$ (cf. also the comment
just after Lemma~\ref{nonzero}).
It is worthwhile noticing that the parameter transformation
defining the dual parameters in \eqref{hattr} is an involution
(the duals of $\hat{t}_r$ bring us back
to $t_r$, $r=0,\ldots,3$) and, furthermore, that
the dual parameters $\hat{t}_r$ satisfy the truncation condition \eqref{trunc}
if the parameters $t_r$ do so (because $\hat{t}_a\hat{t}_b=t_at_b$).
We now form discrete Harish-Chandra-like
$c$-functions $\hat{C}_\pm^{\text{\tiny qR}}$ that
are dual to the $c$-functions
$C_\pm^{\text{\tiny qR}}$ \eqref{Cqr+}, \eqref{Cqr-}
(i.e., $\hat{C}_\pm^{\text{\tiny qR}}$ is obtained from
$C_\pm^{\text{\tiny qR}}$ by replacing the parameters
$t_0,\ldots ,t_3$ by the dual parameters
$\hat{t}_0,\ldots ,\hat{t}_3$)
\begin{eqnarray}\label{hatCqr+}
\hat{C}_+^{\text{\tiny qR}}(\lambda )\!\!\! &=&\!\!\! \hat{c}_0(\lambda)
\prod_{1\leq j< k\leq n}
\left( \frac{(\hat{\tau}_j\hat{\tau}_k;q)_{\lambda_j+\lambda_k}}
            {(t\hat{\tau}_j\hat{\tau}_k;q)_{\lambda_j+\lambda_k}}
\frac{(\hat{\tau}_j\hat{\tau}_k^{-1};q)_{\lambda_j-\lambda_k}}
     {(t\hat{\tau}_j\hat{\tau}_k^{-1};q)_{\lambda_j-\lambda_k}} \right) \\
&& \makebox[2em]{}\times \prod_{1\leq j\leq n}
\left( \frac{(\hat{\tau}_j^2;q)_{2\lambda_j}}
     {\prod_{0\leq r\leq 3}(\hat{t}_r\hat{\tau}_j;q)_{\lambda_j}} \right)
\nonumber\\
\label{hatCqr-}
\hat{C}_-^{\text{\tiny qR}}(\lambda )\!\!\! &=&\!\!\! \hat{c}_0(\lambda)
\prod_{1\leq j< k\leq n}
\left( \frac{(t^{-1}q\hat{\tau}_j\hat{\tau}_k ;q)_{\lambda_j+\lambda_k}}
     {(q\hat{\tau}_j\hat{\tau}_k ;q)_{\lambda_j+\lambda_k}}
\frac{(t^{-1}q\hat{\tau}_j\hat{\tau}_k^{-1} ;q)_{\lambda_j-\lambda_k}}
     {(q\hat{\tau}_j\hat{\tau}_k^{-1} ;q)_{\lambda_j-\lambda_k}} \right) \\
& &\makebox[2em]{} \times \prod_{1\leq j\leq n}
\left( \frac{\prod_{0\leq r\leq 3}(\hat{t}_r^{-1}q\hat{\tau}_j;q)_{\lambda_j}}
     {(q\hat{\tau}_j^2;q)_{2\lambda_j}} \right) \nonumber
\end{eqnarray}
with
\begin{equation}\label{hatc0}
\hat{c}_0(\lambda) = \prod_{1\leq j\leq n} \left( t^{n-j}
(\hat{t}_0\hat{t}_1\hat{t}_2\hat{t}_3 q^{-1})^{1/2} \right)^{\lambda_j} .
\end{equation}
and (cf. \eqref{tau})
\begin{equation}\label{hattau}
\hat{\tau}_j = t^{n-j}\hat{t}_a =
              t^{n-j} (t_0 t_1 t_2 t_3 q^{-1})^{1/2}.
\;\;\;\;\;\;\;\; j=1,\ldots ,n.
\end{equation}
It is important to convince oneself that, despite the
appearances of square roots in the definitions of
$\hat{t}_r$ \eqref{hattr} and $\hat{\tau}_j$ \eqref{hattau},
the functions $\hat{C}_\pm^{\text{\tiny qR}}(\lambda)$
\eqref{hatCqr+}, \eqref{hatCqr-} (including the common factor
$\hat{c}_0$ \eqref{hatc0}) are rational
in the parameters $t$, $t_0,\ldots ,t_3$ and $q$.
The point is that in the above expressions
the quantities $\hat{t}_r$ and $\hat{\tau}_j$ always
occur in rational combinations
of the form $\hat{\tau}_j\hat{\tau}_k^{\pm 1}$,
$\hat{t}_r^{\pm 1}\hat{\tau}_j$, $\hat{\tau_j}^2$
and that $(\hat{t}_0\hat{t}_1\hat{t}_2\hat{t}_3 q^{-1})^{1/2}=
(\hat{t}_a\hat{t}_b\hat{t}_c\hat{t}_d q^{-1})^{1/2}=
t_a$. In this sense the square roots in $\hat{t}_r$
\eqref{hattr} have
merely a formal meaning and these dual parameters were
introduced mostly for notational convenience
and to emphasize the duality between the expressions
for the above $c$-functions
$\hat{C}_\pm^{\text{\tiny qR}}$ and the $c$-functions
$C_\pm^{\text{\tiny qR}}$ appearing in the previous section.

Our main object of interest in this section will be a function
$\mathcal{N}^{\text{\tiny qR}} (\lambda)$ on
the cone $\Lambda$ \eqref{cone} defined by
\begin{equation}\label{NqR}
\mathcal{N}^{\text{\tiny qR}} (\lambda) =
\frac{\hat{C}_-^{\text{\tiny qR}} (\lambda)}
     {\hat{C}_+^{\text{\tiny qR}} (\lambda)}.
\end{equation}
It is clear from (the dual version of)
Lemma~\ref{nonzero} (part a.) and the above comments that
$\hat{C}_\pm^{\text{\tiny qR}}(\lambda )$ and
hence $\mathcal{N}^{\text{\tiny qR}}(\lambda)$ are well-defined
nonzero rational expressions in the parameters $t$, $t_0,\ldots t_3$
and $q$ for all $\lambda\in\Lambda$. Furthermore, since
$\hat{t}_a\hat{t}_b=t_at_b$ the truncation condition \eqref{trunc}
reads the same in the dual parameters $\hat{t}_r$
as in the original parameters $t_r$. Thus, we may apply again (the dual
version) of Lemma~\ref{nonzero} (part b.) to conclude
that for $\lambda\in\Lambda_N$ these expressions are well-defined and
nonzero as rational expressions in the parameters
$t$, $t_0,\ldots t_3$ and $q$ subject to the truncation condition
\eqref{trunc}.

The following theorem provides a summation formula
expressing $\langle p_\lambda ,p_\lambda \rangle_N^{\text{\tiny qR}}$ in terms
of $\langle 1 ,1 \rangle_N^{\text{\tiny qR}}$
(which corresponds to $\lambda= 0$).
\begin{theorem}\label{sumthm}
For parameters subject to the truncation condition
\eqref{trunc} one has that
\begin{equation}\label{sumeq}
\langle p_\lambda , p_\lambda \rangle_N^{\text{\tiny qR}} =
\mathcal{N}^{\text{\tiny qR}}(\lambda)\langle 1, 1\rangle_N^{\text{\tiny qR}}
\;\;\;\;\;\;\;\;\;\;\;\;\; (\lambda\in\Lambda_N),
\end{equation}
with $\mathcal{N}^{\text{\tiny qR}}(\lambda)$ given by \eqref{NqR} and
\begin{equation}\label{unitnorm}
\langle 1, 1\rangle_N^{\text{\tiny qR}} =\sum_{\nu\in\Lambda_N}
\Delta^{\text{\tiny qR}}(\nu)
\end{equation}
(where $\langle \cdot ,\cdot \rangle_N^{\text{\tiny qR}}$ and
$\Delta^{\text{\tiny qR}}(\nu)$
are defined by \eqref{discreteIP} and \eqref{weightqR}, respectively).
\end{theorem}
For an outline of the proof and further details regarding the
Pieri type recurrence formulas lying at the basis of this proof the reader
is referred to Appendix~\ref{appb}.

It is clear that the summation formula \eqref{sumeq}
should be thought of as a set of identities for the multivariable
Askey-Wilson polynomials with parameters
subject to the truncation condition \eqref{trunc}, which
complements the orthogonality identities
in \eqref{orthoqR}.
For $n=1$ the formula reduces to the known
(cf. \cite{ask-wil:set,gas-rah:basic}) summation formula
\begin{equation}
\sum_{0\leq \nu \leq N}
\bigl( p_\lambda (t_a q^{\nu}) \bigr)^2  \Delta^{\text{\tiny qR}} (\nu )=
\frac{(q;q)_\lambda \prod_{0\leq r < s\leq 3} (t_rt_s;q)_\lambda}
     {(t_0t_1t_2t_3q^{\lambda-1};q)_\lambda\; (t_0t_1t_2t_3;q)_{2\lambda}}\;
\langle 1, 1\rangle_N^{\text{\tiny qR}}
\end{equation}
($\lambda = 0,\ldots ,N$), with $\Delta^{\text{\tiny qR}} (\nu )$
taken from \eqref{deltaqR} and
\begin{equation}
\langle 1, 1\rangle_N^{\text{\tiny qR}}=\sum_{0\leq \nu \leq
N}\Delta^{\text{\tiny qR}} (\nu ) ,
\end{equation}
for the monic $q$-Racah polynomial given by
$p_\lambda (z)$ \eqref{bhrepAW}
with parameters subject to the condition
$t_at_b=q^{-N}$.

Although in principle the expression in \eqref{unitnorm} is
in itself already explicit as a finite sum of explicitly
given terms (this in contrast to the sum
$\langle p_\lambda,p_\lambda\rangle_N^{\text{\tiny qR}}$ with $\lambda \neq 0$,
for which in general one does not know the terms in a very
explicit form unless $n=1$),
it would be interesting to also
evaluate the sum $\langle 1,1\rangle_N^{\text{\tiny qR}}$
in terms of a product formula.
For n=1 this was done by Askey and Wilson in \cite{ask-wil:set}
by means of a $q$-Dougall summation formula for a very-well-poised
(terminating) ${}_6\phi_5$ series (see also \cite{gas-rah:basic}),
entailing
\begin{equation}
\sum_{0\leq \nu \leq N}\Delta^{\text{\tiny qR}} (\nu ) =
\frac{(\hat{t}_0\hat{t}_1\hat{t}_2\hat{t}_3,\hat{t}_a^{-1}\hat{t}_b;q)_N}
     {(\hat{t}_b\hat{t}_c,\hat{t}_b\hat{t}_d ;q)_N}=
\frac{(t_a^2q,t_c^{-1}t_d^{-1}q;q)_N}
     {(t_at_c^{-1}q, t_at_d^{-1}q ;q)_N}
\end{equation}
(where $\Delta^{\text{\tiny qR}}(\nu)$ is given by \eqref{deltaqR} and
$t_at_b=q^{-N}$).
For arbitrary number of variables summing
$\langle 1,1\rangle_N^{\text{\tiny qR}}$
amounts to the evaluation of a finite Selberg type $q$-Jackson integral
(related to the root system $BC_n$).
Similar (but infinite) Selberg type $q$-Jackson integrals appear in the recent
works of Aomoto and Ito \cite{aom:product,ito:theta}.
More specifically, the $q$-Jackson integral corresponding
to $\langle 1,1\rangle_N^{\text{\tiny qR}}$ constitutes a finite analogue
of an infinite $q$-Jackson Selberg integral related
to the root system $BC_n$ that belongs to
the same class as those
studied in \cite{aom:product,ito:theta}.
The connection with the $q$-Jackson integrals considered
by Aomoto and Ito becomes particularly transparent after recalling
(see Remark~\ref{cext})
that for $0<|q|<1$ the weights $\Delta^{\text{\tiny qR}}(\nu)$ may be thought
of as the restriction of a function
$\bar{\Delta}^{\text{\tiny qR}}(z)$ of the form in
\eqref{context} to the grid points $\tau q^{\nu}$, $\nu\in\Lambda_N$.

We see from Theorem~\ref{sumthm}
that the quantities
$\langle p_\lambda , p_\lambda\rangle_N^{\text{\tiny qR}}$,
$\lambda\in\Lambda_N$
are nonzero rational expressions in the parameters $t$,
$t_0,\ldots ,t_3$ and $q$ (subject the truncation condition \eqref{trunc})
because the factors
$\mathcal{N}^{\text{\tiny qR}}(\lambda )$ are nonzero.
(Clearly the common factor
$\langle 1 , 1 \rangle_N^{\text{\tiny qR}}$ in \eqref{sumeq}
is nonzero because otherwise the bilinear form
$\langle\cdot ,\cdot \rangle_N^{\text{\tiny qR}}$ would vanish
identically on the space $\mathcal{H}_N^{\text{\tiny qR}}$.)
This checks
with Proposition~\ref{nondegprp}, which stated that
for generic parameters (subject to the truncation condition)
the bilinear form $\langle\cdot ,\cdot \rangle_N^{\text{\tiny qR}}$
is nondegenerate on the space $\mathcal{H}_N^{\text{\tiny qR}}$.

The orthonormalization constants are given by square roots of
$\langle p_{\lambda} , p_\lambda \rangle_N^{\text{\tiny qR}}$,
$\lambda\in\Lambda_N$.
After division by these constants the basis
$p_\lambda$, $\lambda\in\Lambda_N$ becomes an orthonormal
basis for the space $\mathcal{H}_N^{\text{\tiny qR}}$
with respect to
the bilinear form $\langle\cdot ,\cdot \rangle_N^{\text{\tiny qR}}$.
The fact that (generically) the orthogonal basis
$p_\lambda$, $\lambda\in\Lambda_N$ can be turned into
an orthonormal basis for $\mathcal{H}_N^{\text{\tiny qR}}$
(or, equivalently, that the quantities
$\langle p_{\lambda} , p_\lambda \rangle_N^{\text{\tiny qR}}$,
$\lambda\in \Lambda_N$ are nonzero)
implies that (generically) no linear dependences arise between the functions
$p_\lambda$, $\lambda\in\Lambda_N$ when being
restricted to the grid $\tau q^{\nu}$, $\nu\in\Lambda_N$.
This is of course precisely what we already saw
in Lemma~\ref{replem},
and what---because of dimensional considerations---implied
that (generically) any function over the grid could be represented
exactly by a unique function in $\mathcal{H}_N^{\text{\tiny qR}}$
(Proposition~\ref{repprp}).

We also saw already in Section~\ref{sec3} that for parameters subject
to the truncation condition the weights
$\Delta^{\text{\tiny qR}}(\nu)$ \eqref{weightqR}
vanish when $\nu\in \Lambda \setminus \Lambda_N$, because
then the $c$-function $C_+^{\text{\tiny qR}}(\nu)$ \eqref{Cqr+}
became infinite. Since the truncation condition \eqref{trunc}
reads the same in the dual parameters $\hat{t}_r$ \eqref{hattr}
as in the original parameters $t_r$
(recal $\hat{t}_a\hat{t}_b=t_at_b$),
the $c$-function $\hat{C}_+^{\text{\tiny qR}}(\lambda)$ \eqref{hatCqr+}
also becomes infinite for parameters subject to the truncation
condition when $\lambda\in \Lambda \setminus \Lambda_N$. Hence, the sum
$\langle p_\lambda , p_\lambda \rangle_N^{\text{\tiny qR}}$ \eqref{sumeq}
vanishes if the parameters satisfy the truncation condition \eqref{trunc}
when $\lambda\in \Lambda \setminus \Lambda_N$.
The vanishing of the sum
$\langle p_\lambda ,p_\lambda \rangle_N^{\text{\tiny qR}}$
in this situation should
of course not come as a surprise in view of the previously noted fact that
the polynomials in question are actually zero
on the grid points $\tau q^{\nu}$, $\nu\in\Lambda_N$
for these parameters
(Proposition~\ref{findimredprp}).

It is possible to capture the orthogonality relations \eqref{orthoqR}
together with the orthonormalization formulas \eqref{sumeq}
in terms of a purely linear-algebraic formulation.
In this language Theorem~\ref{orthothm} and Theorem~\ref{sumthm}
boil down to the property of the matrix
\begin{equation}\label{qmatrixortho}
[K_{\mu ,\nu}]_{\mu ,\nu\in\Lambda_N}\;\;\;\;\;
\text{with}\;\;\;
K_{\mu ,\nu} \equiv P_{\mu} (\tau q^{\nu})
\frac{\Bigl( \hat{\Delta}^{\text{\tiny qR}}(\mu) \Bigr)^{1/2}
      \Bigl( \Delta^{\text{\tiny qR}}(\nu) \Bigr)^{1/2}}
     {\Bigl( \langle 1 , 1 \rangle_N^{\text{\tiny qR}}\Bigr)^{1/2} }
\end{equation}
being an orthogonal matrix for generic parameters subject
to the truncation condition \eqref{trunc}.
Here we have used the renormalized polynomial
$P_{\mu} (z)\equiv \hat{C}_+^{\text{\tiny qR}}(\mu)p_\mu (z)$
(cf. Appendix~\ref{appb}) and the weights for the
`Plancherel' measure (cf. \eqref{weightqR})
\begin{equation*}
\hat{\Delta}^{\text{\tiny qR}}(\mu)=
\frac{1}{\hat{C}_+^{\text{\tiny qR}}(\mu)\hat{C}_-^{\text{\tiny qR}}(\mu)}
\end{equation*}
and, we have again (cf. \eqref{qmatrix})
ordered the rows
and columns of the matrix according to an arbitrary total extension
of the partial order in \eqref{order}.
The orthogonality (and thus invertibility) of the matrix
$K=[K_{\mu ,\nu}]_{\mu ,\nu\in\Lambda_N}$ in \eqref{qmatrixortho}
checks with the previously noted invertibility of
the matrix $M$ given by \eqref{qmatrix} because---as should be clear from
our discussion--- one has that $K=LMD$ with $L$ triangular and
$D$ diagonal and both $L$ and $D$ (generically) nonsingular (as the
matrix elements on their respective
diagonals are nonzero as meromorphic functions in the parameters
subject to the truncation condition).

\begin{remark}\label{dualrem}
The normalization of the polynomials
$P_\mu (z)=\hat{C}_+^{\text{\tiny qR}}(\mu)p_\mu (z)$ is such
that $P_\mu (\tau )=1$ (see Appendix~\ref{appb}).
Moreover, the
renormalized multivariable Askey-Wilson polynomials $P_\mu (z)$ satisfy
Macdonald's duality property (see also Appendix~\ref{appb}) stating that
\begin{equation}\label{dualprop}
P_{\mu} (\tau q^{\nu})=\hat{P}_{\nu} (\hat{\tau} q^{\mu}),
\end{equation}
where $\hat{P}_{\nu} (z)=C_+^{\text{\tiny qR}}(\nu)\hat{p}_\nu (z)$
denotes the renormalized multivariable Askey-Wilson polynomial
dual to $P_{\nu}(z)=\hat{C}_+^{\text{\tiny qR}}(\nu)p_\nu (z)$, i.e.
with the parameters $t_r$ being replaced by the dual parameters
$\hat{t}_r$ \eqref{hattr}. (For $\nu =0$ the equality in \eqref{dualprop}
reduces to the evaluation formula $P_\mu (\tau )=1$.)
The duality relation \eqref{dualprop} for the polynomials
gives rise to a similar duality property for the matrix
$[K_{\mu ,\nu}]_{\mu ,\nu\in\Lambda_N}$ \eqref{qmatrixortho},
viz., transposition of $[K_{\mu ,\nu}]_{\mu ,\nu\in\Lambda_N}$
leads one to the matrix $[\hat{K}_{\mu ,\nu}]_{\mu ,\nu\in\Lambda_N}$
in which the parameters $t_r$ are replaced by the dual
parameters $\hat{t}_r$. Here we also used that the orthogonality
of the matrices fixes their normalization uniquely.
In particular one sees that as a consequence the sum
$\langle 1 ,1 \rangle^{\text{\tiny qR}}_N$ is invariant with
respect to the transformation $t_r\rightarrow \hat{t}_r$ (for
parameters satisfying the truncation condition \eqref{trunc}).
For $n=1$ this property of $\langle 1 ,1 \rangle^{\text{\tiny qR}}_N$
boils down to the equality
\begin{equation}
\frac{(\hat{t}_a^2q,\hat{t}_c^{-1}\hat{t}_d^{-1}q;q)_N}
     {(\hat{t}_a\hat{t}_c^{-1}q, \hat{t}_a\hat{t}_d^{-1}q ;q)_N}=
\frac{(t_a^2q,t_c^{-1}t_d^{-1}q;q)_N}
     {(t_at_c^{-1}q, t_at_d^{-1}q ;q)_N},
\end{equation}
which is not difficult to check directly
using  \eqref{trunc}, \eqref{hattr}
and the transformation property for $q$-shifted factorials
\begin{equation*}
(a_1q^{-N};q)_N/(a_2q^{-N};q)_N =
(a_1/a_2)^N (a_1^{-1}q;q)_N/ (a_2^{-1}q;q)_N.
\end{equation*}
In the special case that the parameters satisfy the additional constraint
$t_aq=t_bt_ct_d$ one has that $\hat{t}_r=t_r$ ($r=0,\ldots ,3$).
Hence, we are then in a self-dual situation
$\hat{P}_{\mu}(z)=P_{\mu}(z)$ and the matrix
$[K_{\mu ,\nu}]_{\mu ,\nu\in\Lambda_N}$
\eqref{qmatrixortho}---in addition to being orthogonal---now also
becomes symmetric.
\end{remark}

\section{Transition to Racah type polynomials}\label{sec5}
We will now study the transition from the {\em basic} hypergeometric
level to the hypergeometric level.
To this end we substitute the variables
\begin{equation}\label{varsub}
z_j=q^{x_j},\;\;\;\;\;\;\;\;\; j=1,\ldots ,n
\end{equation}
and in addition perform a reparametrization of the form
\begin{equation}\label{repa}
t=q^g,\;\;\;\;\;\;\; t_r = q^{g_r},\;\;\; r=0,\ldots ,3.
\end{equation}
After these substitutions and division by a constant factor
$(1-q)^2$, the $q$-difference operator $D$ \eqref{DAW}
passes for $q\rightarrow 1$
over into a second order difference operator given by
\begin{equation}\label{DW}
\tilde{D} = \sum_{1\leq j\leq n} \Bigl(
\tilde{V}_j(x) (T_{j}-1) \; +\; \tilde{V}_{-j}(x) (T_{j}^{-1}-1) \Bigr)
\end{equation}
where
\begin{eqnarray*}
\tilde{V}_{\varepsilon j}(x) &=&
\frac{ \prod_{0\leq r\leq 3} (g_r+\varepsilon x_j)}
     { (2\varepsilon x_j )\, (1 +2\varepsilon x_j )} \\
& & \times \prod_{\stackrel{1\leq k \leq n}{ k\neq j}}
\left( \frac{g +\varepsilon x_j +x_k}{\varepsilon x_j+x_k}\right)
\left( \frac{g +\varepsilon x_j -x_k}{\varepsilon x_j-x_k}\right) ,
\;\;\;\;\;\;\;\varepsilon =\pm 1
\end{eqnarray*}
and with the action of the operators $T_j$, $j=1,\ldots ,n$ being of the form
\begin{equation*}
(T_{j}f) (x_1,\ldots ,x_n) =
f(x_1,\ldots ,x_{j-1},x_j+1,
x_{j+1}, \ldots ,x_n).
\end{equation*}
It turns out (see \cite{die:properties})
that the difference operator $\tilde{D}$ \eqref{DW} is triangular
with respect to the (partially ordered) basis of symmetrized
monomials for the space
$\mathbb{C}^{S_n}[x_1^2,\ldots ,x_n^2]$
(consisting of the
permutation-invariant and even polynomials in the
variables $x_1,\ldots ,x_n$). Specifically, one has that
\begin{equation}\label{triaW}
\tilde{D} \tilde{m}_{\lambda} = \sum_{\mu\in\Lambda,\;\mu\leq\lambda}
\tilde{E}_{\lambda,\mu}\, \tilde{m}_{\mu}\;\;\;\;\text{with}\;\;\;\;
\tilde{E}_{\lambda,\mu}\in \mathbb{C}[g,g_0,g_1,g_2,g_3]
\end{equation}
where
\begin{equation}\label{monomW}
\tilde{m}_\lambda (x) =
\sum_{\mu\in S_n (\lambda )} x_1^{2\mu_1}\cdots x_n^{2\mu_n},
\;\;\;\;\;\;\;\;\;\; \lambda\in\Lambda .
\end{equation}
Here the summation in \eqref{monomW} is meant over the orbit
of $\lambda\in\Lambda$ \eqref{cone}
under the action of the permutation group $S_n$
(which permutes the vector components $\lambda_1,\ldots ,\lambda_n$)
and the partial order of the cone $\Lambda$ \eqref{cone} is
taken to be the same as before (see \eqref{order}).
The diagonal matrix elements $\tilde{E}_{\lambda ,\lambda}$
in \eqref{triaW}
(which can be obtained for
$q\rightarrow 1$ from $E_{\lambda ,\lambda}$
\eqref{eigenvAW} after substitution of
\eqref{repa} and division by $(1-q)^2$)
read explicitly
\begin{equation}\label{eigenvW}
\tilde{E}_{\lambda ,\lambda} =
\sum_{1\leq j\leq n}
\bigr( (\lambda_j+\hat{\rho}_j)^2-\hat{\rho}_j^2 \bigl) ,
\;\;\;\;\;\;\;\;\; \lambda\in\Lambda
\end{equation}
with
\begin{equation}\label{hatrho}
\hat{\rho}_j = (n-j)g + (g_0+g_1+g_2 +g_3-1)/2,\;\;\;\;\;\; j=1,\ldots ,n.
\end{equation}
The triangularity of the difference operator again reduces the corresponding
eigenvalue problem in the space of the
permutation-invariant and even polynomials
to an in essence finite-dimensional problem.
Although the eigenvalues \eqref{eigenvW} are no
longer nondegenerate, it still remains true that
$\tilde{E}_{\lambda ,\lambda} \neq \tilde{E}_{\mu ,\mu}$
as polynomial expression in the parameters
$g$, $g_0,\ldots ,g_3$, if $\lambda\neq \mu$ and $\lambda, \mu$ are comparable
with respect to the partial order \eqref{order}.
Fortunately, this is already
sufficient to single out the eigenfunctions uniquely
by means of conditions analogous to those
entering the definition of the multivariable Askey-Wilson
polynomials $p_\lambda$ in Section~\ref{sec2}.
It turns out (see \cite{die:properties})
that in the present case we are in fact dealing with
a multivariable
analogue of the Wilson polynomials \cite{wil:some}.
\begin{definition}
The {\em multivariable Wilson polynomial}
associated with a (dominant weight) vector $\lambda\in \Lambda$ \eqref{cone}
is the (unique) monic permutation
invariant and even polynomial of the form
\begin{subequations}
\begin{equation}\label{con1}
\tilde{p}_{\lambda}(x)=\tilde{m}_{\lambda}(x)+
       \sum_{\mu\in\Lambda ,\;\mu<\lambda}c_{\lambda,\mu}\tilde{m}_{\mu}(x)
\;\; \text{with}\;\; c_{\lambda,\mu}\in {\mathbb C}(g,g_0,g_1,g_2 ,g_3),
\end{equation}
such that
\begin{equation}\label{con2}
\tilde{D}\, \tilde{p}_{\lambda}=\tilde{E}_{\lambda,\lambda}\,
\tilde{p}_{\lambda} .
\end{equation}
\end{subequations}
\end{definition}

We can again represent these polynomials
in terms of a formula of the type \eqref{Macfor}:
\begin{equation}\label{MacforW}
\tilde{p}_{\lambda} = \Biggl( \;\prod_{\mu\in\Lambda,\; \mu<\lambda}
\frac{\tilde{D}-\tilde{E}_{\mu,\mu}}
     {\tilde{E}_{\lambda,\lambda}-\tilde{E}_{\mu,\mu}}\Biggr)\;
\tilde{m}_{\lambda}  .
\end{equation}
Furthermore, one may always replace the monomial basis
in the r.h.s. of such a formula
by any other basis related to it via a unitriangular transformation
(this is immediate from the comments following \eqref{Macfor}
that proved the validity of this type of representations
for the polynomials).
In \cite{die:properties} it was demonstrated, by performing
a suitable unitriangular transformation of the basis elements, that
formula \eqref{Macfor} for the multivariable Askey-Wilson polynomials
tends to formula \eqref{MacforW} for the multivariable
Wilson polynomials in the limit $q\rightarrow 1$ after
substitution of \eqref{varsub}, \eqref{repa} and division by
a constant factor $(1-q)^{2 |\lambda |}$. Here we have
used the (standard) notation
\begin{equation}
|\lambda | \equiv \lambda_1+\cdots +\lambda_n.
\end{equation}
(The point of the unitriangular transformation is that the
basis elements \eqref{monom} all collapse into
constant functions for $q\rightarrow 1$ after
the substitution \eqref{varsub};
by temporarily passing to a basis of $\mathcal{H}^W$
with elements of the form
\begin{equation}\label{modbasis}
\sum_{\mu\in S_n(\lambda )}
(z_1+z_1^{-1}-2)^{\mu_1}\cdots (z_n+z_n^{-1}-2)^{\mu_n},
 \;\;\;\;\;\;\;\;\lambda\in \Lambda ,
\end{equation}
it is seen that a nontrivial limit is obtained after one divides
out the constant factor $(1-q)^{2|\lambda |}$.)

The upshot is that there exists the following limiting relation
between
the multivariable Askey-Wilson polynomials $p_{\lambda}$
of Section~\ref{sec2}
and their Wilson type counterparts $\tilde{p}_\lambda$ of the
present section \cite{die:properties}.
\begin{proposition}\label{limprp}
For Askey-Wilson parameters given by \eqref{repa}, one has that
\begin{equation}
\tilde{p}_\lambda (x) =  \lim_{q\rightarrow 1}
(1-q)^{-2|\lambda |}\; p_\lambda(q^x)
\;\;\;\;\;\;\;\;\;\;\;\;\;\;\;\;
\lambda\in\Lambda
\end{equation}
(where $q^x\equiv (q^{x_1},\ldots ,q^{x_n})$).
\end{proposition}

In \cite{die:properties} the orthogonality
properties of polynomials $\tilde{p}_\lambda (x)$ were investigated
with respect to a continuous Wilson type weight function
$\Delta^{\text{\tiny W}}$. For parameters satisfying
\begin{equation}\label{Wcond}
g \geq 0, \;\;\;\;\;\; \text{Re}(g_r) > 0  \;\; (r=0,1,2 ,3),
\end{equation}
with possible non-real parameters $g_r$ occurring in complex conjugate pairs,
the relevant orthogonality relations read
\begin{equation}
\int_{-\infty}^{\infty}\!\!\!\!\!\!\cdots
               \int_{-\infty}^{\infty}
\tilde{p}_\lambda (ix)\, \overline{\tilde{p}_{\mu}(ix)}\, \Delta^W (x)\,
dx_1\cdots dx_n =0\;\;\; \text{if}\;\;\; \lambda \neq \mu
\end{equation}
where
\begin{eqnarray*}
\Delta^{\text{\tiny W}} (x) &=&
\prod_{\stackrel{1\leq j< k \leq n}
                {\varepsilon_1 ,\varepsilon_2 =\pm 1}}
\frac{\Gamma (g + i(\varepsilon_1 x_j
           +\varepsilon_2 x_{k}) )}
     {\Gamma (i(\varepsilon_1 x_j
           +\varepsilon_2 x_{k}) )} \\
& &  \times \prod_{\stackrel{1\leq j\leq n}{\varepsilon =\pm 1}}
\frac{\Gamma (g_0+i\varepsilon x_j)
\Gamma (g_1+i\varepsilon x_j)
\Gamma (g_2+i\varepsilon x_j)
\Gamma (g_3+i\varepsilon x_j)}{\Gamma (2i\varepsilon x_j )} \nonumber
\end{eqnarray*}
(with $\Gamma (\cdot)$ denoting the gamma function).

We will now apply the limiting relation of Proposition~\ref{limprp}
to the results of Section~\ref{sec3} and \ref{sec4} to
infer that for generic parameters subject to the truncations condition
\begin{equation}\label{truncdeg}
(n-1)g + g_a +g_b + N =0
\end{equation}
(with $N$ a nonnegative integer and
$a,b\in \{ 0,1,2,3\}$ such that $a\neq b$)
the polynomials $\tilde{p}_\lambda (x)$, $\lambda\in\Lambda_N$ \eqref{alcove}
constitute an orthogonal basis for the finite-dimensional space
\begin{equation}
\mathcal{H}^{\text{\tiny R}}_N \equiv
\text{Span} \{ \tilde{m}_\lambda \}_{\lambda \in\Lambda_N}
\end{equation}
endowed with a nondegenerate bilinear form determined by
\begin{equation}\label{discreteIPdeg}
\langle f ,g \rangle_N^{\text{\tiny R}} =
\sum_{\nu \in\Lambda_N} f (\rho+\nu)
 g (\rho+\nu) \Delta^{\text{\tiny R}} (\nu)
\;\;\;\;\;\;\;\;\;\;\; (f,g\in\mathcal{H}^{\text{\tiny R}}_N).
\end{equation}
Here the vector $\rho$ is of the form
(c.f. also the vector $\rho$ in Remark~\ref{gridstruct}
with $\tau$ given by \eqref{tau} and parameters taken from \eqref{repa})
\begin{equation}\label{rho}
\rho = (\rho_1,\ldots ,\rho_n)\;\;\; \text{with}\;\;\;
\rho_j = (n-j)g +g_a,\;\;\; j=1,\ldots ,n
\end{equation}
and the weights are given by
\begin{equation}\label{weightR}
\Delta^{\text{\tiny R}} (\nu ) =
\frac{ 1 }
{C_+^{\text{\tiny R}}(\nu )\: C_-^{\text{\tiny R}}(\nu ) } ,
\end{equation}
with
\begin{eqnarray}\label{Cr+}
C_+^{\text{\tiny R}}(\nu )\!\!\! &=&\!\!\!
\prod_{1\leq j< k\leq n}
\left( \frac{(\rho_j+\rho_k)_{\nu_j+\nu_k}}
            {(g+\rho_j+\rho_k)_{\nu_j+\nu_k}}
\frac{(\rho_j-\rho_k)_{\nu_j-\nu_k}}
     {(g+\rho_j-\rho_k)_{\nu_j-\nu_k}} \right) \\
&& \makebox[2em]{}\times \prod_{1\leq j\leq n}
\left( \frac{(2\rho_j)_{2\nu_j}}
     {\prod_{0\leq r\leq 3}(g_r+\rho_j)_{\nu_j}} \right) ,\nonumber\\
\label{Cr-}
C_-^{\text{\tiny R}}(\nu )\!\!\! &=&\!\!\!
\prod_{1\leq j< k\leq n}
\left( \frac{(1-g+\rho_j+\rho_k )_{\nu_j+\nu_k}}
            {(1+\rho_j+\rho_k )_{\nu_j+\nu_k}}
       \frac{(1-g+\rho_j-\rho_k )_{\nu_j-\nu_k}}
            {(1+\rho_j-\rho_k )_{\nu_j-\nu_k}} \right) \\
& &\makebox[2em]{} \times \prod_{1\leq j\leq n}
\left( \frac{\prod_{0\leq r\leq 3}(1-g_r+\rho_j)_{\nu_j}}
            {(1+2\rho_j)_{2\nu_j}} \right) ,\nonumber
\end{eqnarray}
where we have used Pochhammer symbols defined by
\begin{eqnarray*}
(a_1,\ldots ,a_s)_m &=& (a_1)_m\cdots (a_s)_m, \\
(a)_m &=& a(a+1)\cdots (a+m-1)
\end{eqnarray*}
(with $(a)_0\equiv 1$).

In order to describe the corresponding orthonormalization constants
we again need dual parameters
\begin{equation}\label{hatgr}
\left(
\begin{matrix} \hat{g}_a \\ \hat{g}_b \\ \hat{g}_c \\ \hat{g}_d \end{matrix}
                                                              \right)
= \frac{1}{2} \left(
\begin{array}{rrrr}
1  &  1  &  1  &  1  \\
1  &  1  &  -1  &  -1  \\
1  &  -1  &  1  &  -1  \\
1  &  -1  &  -1  &  1
\end{array}             \right)
\left( \begin{matrix} g_a \\ g_b \\ g_c \\ g_d \end{matrix} \right)
+\frac{1}{2}
\left( \begin{matrix} -1 \\ 1 \\ 1 \\ 1 \end{matrix} \right)
\end{equation}
($\{ a,b,c,d\} =\{ 0,1,2,3 \}$)
and a function $\mathcal{N}^{\text{\tiny R}}(\lambda )$
on the cone $\Lambda$ \eqref{cone}
\begin{equation}\label{NR}
\mathcal{N}^{\text{\tiny R}} (\lambda) =
\frac{\hat{C}_-^{\text{\tiny R}} (\lambda)}
     {\hat{C}_+^{\text{\tiny R}} (\lambda)}
\end{equation}
that is governed by $c$-functions dual to $C_\pm^{\text{\tiny R}}$
\eqref{Cr+}, \eqref{Cr-}
\begin{eqnarray}\label{hatCr+}
\hat{C}_+^{\text{\tiny R}}(\lambda )\!\!\! &=&\!\!\!
\prod_{1\leq j< k\leq n}
\left( \frac{(\hat{\rho}_j+\hat{\rho}_k)_{\lambda_j+\lambda_k}}
            {(g+\hat{\rho}_j+\hat{\rho}_k)_{\lambda_j+\lambda_k}}
\frac{(\hat{\rho}_j-\hat{\rho}_k)_{\lambda_j-\lambda_k}}
     {(g+\hat{\rho}_j-\hat{\rho}_k)_{\lambda_j-\lambda_k}} \right) \\
&& \makebox[2em]{}\times \prod_{1\leq j\leq n}
\left( \frac{(2\hat{\rho}_j)_{2\lambda_j}}
     {\prod_{0\leq r\leq 3}(\hat{g}_r+\hat{\rho}_j)_{\lambda_j}} \right) ,
\nonumber\\
\label{hatCr-}
\hat{C}_-^{\text{\tiny R}}(\lambda )\!\!\! &=&\!\!\!
\prod_{1\leq j< k\leq n}
\left( \frac{(1-g+\hat{\rho}_j+\hat{\rho}_k )_{\lambda_j+\lambda_k}}
            {(1+\hat{\rho}_j+\hat{\rho}_k )_{\lambda_j+\lambda_k}}
       \frac{(1-g+\hat{\rho}_j-\hat{\rho}_k )_{\lambda_j-\lambda_k}}
            {(1+\hat{\rho}_j-\hat{\rho}_k )_{\lambda_j-\lambda_k}} \right) \\
& &\makebox[2em]{} \times \prod_{1\leq j\leq n}
\left( \frac{\prod_{0\leq r\leq 3}(1-\hat{g}_r+\hat{\rho}_j)_{\lambda_j}}
            {(1+2\hat{\rho}_j)_{2\lambda_j}} \right) .\nonumber
\end{eqnarray}
(Recall that the components
of the vector $\hat{\rho}$ are given by \eqref{hatrho},
so $\hat{\rho}_j = (n-j)g+\hat{g}_a$).

The following theorem, which describes the orthogonality properties
of the polynomials $\tilde{p}_\lambda$, $\lambda\in\Lambda_N$
with respect to the bilinear form
$\langle\cdot ,\cdot\rangle_N^{\text{\tiny R}}$, is an immediate
consequence of the application of Proposition~\ref{limprp}
to Theorem~\ref{orthothm} and Theorem~\ref{sumthm}.
\begin{theorem}\label{orthosumdegthm}
For parameters subject to the truncation condition
\eqref{truncdeg} one has that
\begin{equation}\label{orthodeg}
\langle \tilde{p}_{\lambda} , \tilde{p}_{\mu} \rangle_N^{\text{\tiny R}}= 0
\;\;\; \text{for}\;\;\; \lambda \neq \mu
\;\;\;\;\;\;\;\;\;\;\;\;\;\;\;\;\;\; (\lambda,\mu \in\Lambda_N)
\end{equation}
and that
\begin{equation}\label{sumdeg}
\langle \tilde{p}_{\lambda} , \tilde{p}_{\lambda} \rangle_N^{\text{\tiny R}}
= \mathcal{N}^{\text{\tiny R}}(\lambda)\:
  \langle 1, 1\rangle_N^{\text{\tiny R}}
     \;\;\;\;\;\;\;\;\;\;\;\;\;\;\;\;\; (\lambda\in\Lambda_N)
\end{equation}
with $\mathcal{N}^{\text{\tiny R}}(\lambda)$ given by \eqref{NR} and
\begin{equation}
\langle 1, 1\rangle_N^{\text{\tiny R}} =\sum_{\nu\in\Lambda_N}
\Delta^{\text{\tiny R}}(\nu)
\end{equation}
(where $\langle \cdot ,\cdot \rangle_N^{\text{\tiny R}}$ and
$\Delta^{\text{\tiny R}}(\nu)$
are defined by \eqref{discreteIPdeg} and \eqref{weightR}, respectively).
\end{theorem}
(The formulas \eqref{orthodeg} and \eqref{sumdeg} should again be
interpreted as equalities between expressions that are rational in
parameters subject to the truncation condition.)
To verify the theorem it suffices to infer that for $q\rightarrow 1$
and parameters given
by \eqref{repa} the $c$-functions
$C^{\text{\tiny qR}}_+(\nu )$ \eqref{Cqr+} and
$\hat{C}^{\text{\tiny qR}}_+(\lambda)$ \eqref{hatCqr+}
(multiplied by a factor $(1-q)^{2|\nu |}$ and $(1-q)^{2|\lambda |}$,
respectively)
converge to $C^{\text{\tiny R}}_+(\nu )$ \eqref{Cr+} and
$\hat{C}^{\text{\tiny R}}_+(\lambda)$ \eqref{hatCr+};
and that, similarly,
the $c$-functions
$C^{\text{\tiny qR}}_-(\nu)$ \eqref{Cqr-} and
$\hat{C}^{\text{\tiny qR}}_-(\lambda)$ \eqref{hatCqr-}
(divided by a factor $(1-q)^{2|\nu |}$ and $(1-q)^{2|\lambda |}$,
respectively) tend to
$C^{\text{\tiny R}}_-(\nu)$ \eqref{Cr-} and
$\hat{C}^{\text{\tiny R}}_-(\lambda)$ \eqref{hatCr-} in this limit.
To this end one simply
uses that the (renormalized) $q$-shifted factorial
$(a;q)_m/(1-q)^m$ converges to the Pochhammer symbol $(a)_m$
when $q$ tends to one.
Furthermore, for the parameters \eqref{repa}
the truncation condition \eqref{trunc}
amounts to the condition $(n-1)g +g_a+g_b+N = 0$ (mod
$2\pi i /\log (q)$), which entails \eqref{truncdeg}
in the limit $q\rightarrow 1$.

It is important to again convince oneself that the formulas
\eqref{orthodeg} and
\eqref{sumdeg} are indeed well-defined as rational expressions
in the parameters $g$ and $g_0,\ldots ,g_3$ subject to the truncation
condition \eqref{truncdeg} (i.e., no denominator becomes identical
to zero) and, furthermore, that the r.h.s. of \eqref{sumdeg} is nonzero
as a rational expression in these parameters.
We thus have that the bilinear form
$\langle \cdot ,\cdot \rangle_N^{\text{\tiny R}}$ is
nondegenerate on the space $\mathcal{H}_N^{\text{\tiny R}}$
for generic parameters subject to the truncation condition \eqref{truncdeg}
and that any function defined
on the grid points $\rho +\lambda$, $\lambda\in\Lambda_N$
can be represented exactly by the restriction
to the grid of a unique function in the space
$\mathcal{H}_N^{\text{\tiny R}}$. (The restriction to the grid points
of the orthonormalized basis
$( \langle \tilde{p}_\lambda ,\tilde{p}_\lambda
\rangle_N^{\text{\tiny R}})^{-1/2} \tilde{p}_\lambda (x)$,
$\lambda\in\Lambda_N$
for $\mathcal{H}_N^{\text{\tiny R}}$ yields
an orthonormal
basis for the space of functions over the grid $\rho +\Lambda_N$
endowed with the nondegenerate
bilinear form $\sum_{\nu\in\Lambda_N} f(\rho +\nu) g(\rho +\nu)
\Delta^{\text{\tiny R}}(\nu )$.)
Finally, as degeneration of Proposition~\ref{findimredprp} we arrive at
a similar statement for the polynomials $\tilde{p}_\lambda$.
\begin{proposition}
For parameters subject to the truncation condition
\eqref{truncdeg} one has that
\begin{equation}
\tilde{p}_{\lambda} (\rho+\nu)=0\;\;\;\;\;
\text{for}\;\; \nu\in\Lambda_N
\end{equation}
if $\lambda\in\Lambda\setminus\Lambda_N$.
\end{proposition}

The corresponding linear-algebraic formulation of Theorem~\ref{orthosumdegthm}
states that the matrix (cf. \eqref{qmatrixortho})
\begin{equation}\label{matrixortho}
[\tilde{K}_{\mu ,\nu}]_{\mu ,\nu\in\Lambda_N}\;\;\;\;\;
\text{with}\;\;\;
\tilde{K}_{\mu ,\nu} \equiv \tilde{P}_{\mu} (\tau q^{\nu})
\frac{\Bigl( \hat{\Delta}^{\text{\tiny R}}(\mu) \Bigr)^{1/2}
      \Bigl( \Delta^{\text{\tiny R}}(\nu) \Bigr)^{1/2}  }
     {\Bigl( \langle 1 , 1 \rangle_N^{\text{\tiny R}} \Bigr)^{1/2} }
\end{equation}
is orthogonal for generic parameters $g$, $g_0,\ldots ,g_3$
subject to the truncation condition \eqref{truncdeg}.
Here $\tilde{P}_{\mu}$ is the renormalized polynomial
$\tilde{P}_{\mu} (x)\equiv\hat{C}^{\text{\tiny R}}_+(\mu)\: \tilde{p}_\mu (x)$
and  $\hat{\Delta}^{\text{\tiny R}}$ denotes the
weight function for the discrete `Plancherel' measure
$\hat{\Delta}^{\text{\tiny R}}(\mu )=
1/(\hat{C}^{\text{\tiny R}}_+(\mu)\:\hat{C}^{\text{\tiny R}}_-(\mu) )$
dual to the weight function $\Delta^{\text{\tiny R}}$ \eqref{weightR}.

As a side remark we mention that the orthogonality
of $[\tilde{K}_{\mu ,\nu}]_{\mu ,\nu\in\Lambda_N}$
implies that the matrix
\begin{equation}\label{matrix}
\tilde{M}\equiv [ \tilde{m}_{\mu}(\rho +\nu )]_{\mu ,\nu\in\Lambda_N},
\end{equation}
has a determinant that does not vanish as a polynomial in the
parameters $g$ and $g_a$.
This is because $\tilde{M}$ is related to
$\tilde{K}\equiv [\tilde{K}_{\mu ,\nu}]_{\mu ,\nu\in\Lambda_N}$
\eqref{matrixortho} by $\tilde{K}=\tilde{L}\tilde{M}\tilde{D}$, where
$\tilde{L}$ is triangular and $\tilde{D}$ is diagonal
with both matrices having elements on the diagonal
that do not vanish as meromorphic functions in the
parameters $g$ and $g_0,\ldots ,g_3$.
It is interesting to observe that
the invertibility of the matrix $\tilde{M}$
does not seem so easily established directly
(i.e. without using $\tilde{K}$),
as was the case when dealing with its $q$-version
$M= [m_{\mu}(\tau q^{\nu} )]_{\mu ,\nu\in\Lambda_N}$ \eqref{qmatrix}
just after Lemma~\ref{replem}. (The problem is of course that here
in the degenerate version we have lost the possibility to play
with the parameter $q$.) Notice, however, that in the special case of
only one single variable the matrix $\tilde{M}$ \eqref{matrix} becomes
a Vandermonde matrix $[(g_a +\nu)^{2\mu}]_{0\leq \mu ,\nu \leq N}$,
from which it immediately follows that the determinant is nonzero
for generic parameter values.

The difference equation \eqref{con2} tells us that
for $n=1$ the polynomials $\tilde{p}_\lambda$
reduce to the  monic Wilson polynomials \cite{wil:some}
\begin{eqnarray}\label{hrepW}
\tilde{p}_\lambda (x) &=&
\frac{(g_0+g_1, g_0+g_2, g_0+g_3)_\lambda}
     { (g_0+g_1+g_2+g_3+\lambda -1)_\lambda } \times\\
& & {}_4F_3
\left(
\begin{matrix}
-\lambda ,\; g_0+g_1+g_2+g_3+\lambda -1,\;
g_0+x,\; g_0-x \\ [0.5ex]
g_0+g_1 ,\; g_0+g_2 ,\; g_0+g_3
\end{matrix} \; ; 1\right) , \nonumber
\end{eqnarray}
where we have used standard notation for
the hypergeometric series (see e.g. \cite{ask-wil:some,gas-rah:basic})
\begin{equation*}
{}_rF_s\left( \begin{array}{c}
             a_1,\ldots ,a_r \\ b_1,\ldots ,b_s
            \end{array} ; z \right) =
\sum_{k=0}^\infty
\frac{(a_1,\ldots ,a_r)_k}{(b_1,\ldots ,b_s)_k} \frac{z^k}{k!} .
\end{equation*}
(The explicit hypergeometric
representation \eqref{hrepW} for $\tilde{p}_\lambda$
in the case of one variable
also follows from Proposition~\ref{limprp} and the corresponding basic
hypergeometric formula for $p_\lambda$ in \eqref{bhrepAW}.)
The identities in Theorem~\ref{orthosumdegthm}
amount in this special case to the discrete
orthogonality relations
for the monic Wilson polynomials \eqref{hrepW} subject to the
parameter condition $g_a+g_b+N=0$ \cite{wil:some}
\begin{equation*}
\sum_{0\leq \nu \leq N}
\tilde{p}_\lambda (g_a +\nu) \tilde{p}_{\mu} (g_a +\nu)
\Delta^{\text{\tiny R}} (\nu ) = 0
\;\;\;\text{for} \;\;\;\lambda \neq \mu
\end{equation*}
($\lambda ,\mu \in \{ 0,\ldots ,N\} $) and
\begin{eqnarray*}
&&\sum_{0\leq \nu \leq N}
\tilde{p}_\lambda (g_a +\nu)\, \tilde{p}_{\lambda} (g_a +\nu)\,
\Delta^{\text{\tiny R}} (\nu )\: = \\
&&\makebox[5em]{}
\frac{\lambda ! \prod_{0\leq r< s\leq 3}(g_r+g_s)_\lambda}
     {(g_0+g_1+g_2+g_3+\lambda-1)_\lambda\, (g_0+g_1+g_2+g_3)_{2\lambda}}
\: \langle 1, 1\rangle_N^{\text{\tiny R}}
\end{eqnarray*}
($\lambda \in \{ 0,\ldots ,N\} $), where
\begin{eqnarray*}
\langle 1, 1\rangle_N^{\text{\tiny R}} &=&
         \sum_{0\leq \nu \leq N} \Delta^{\text{\tiny R}} (\nu) ,\\
\Delta^{\text{\tiny R}} (\nu) &=&
                 \left( 1+ \frac{\nu}
                                {g_a}   \right)
\ \frac{(g_0+g_a,g_1+g_a,g_2+g_a,g_3+g_a)_{\nu}}
       {(1-g_0+g_a,1-g_1+g_a,1-g_2+g_a,1-g_3+g_a)_{\nu}} .
\end{eqnarray*}
In the finite-dimensional case with the truncation
condition $g_a+g_b+N=0$ and discrete orthogonality properties,
the one-variable polynomials $\tilde{p}_\lambda$ in
\eqref{hrepW} are usually referred to
as Racah polynomials rather than Wilson polynomials.
The normalization factor $\langle 1, 1\rangle_N^{\text{\tiny R}}$
for the one-variable Racah polynomials can be evaluated
in product form by means of
a summation formula for a very-well-poised
(terminating) ${}_5F_4$ series
due to Dougall, which entails \cite{wil:some} (see also
\cite{gas-rah:basic} for the Dougall ${}_5F_4$ summation formula)
\begin{equation}
\sum_{0\leq \nu \leq N} \Delta^{\text{\tiny R}} (\nu) =
\frac{(1+2g_a,1-g_c-g_d)_N}{(1+g_a-g_c,1+g_a-g_d)_N} .
\end{equation}

\begin{remark}
The duality relations  for the (renormalized)
multivariable Askey-Wils\-on polynomials
in Remark~\ref{dualrem} give in the limit $q\rightarrow 1$ rise to
analogous duality relations for the (renormalized) multivariable
Wilson polynomials
$\tilde{P}_\mu (x)= \hat{C}_+^{\text{\tiny R}}(\mu ) \tilde{p}_\mu (x)$
\cite{die:properties}
\begin{equation}\label{dualpropdeg}
\tilde{P}_{\mu} (\rho + \nu )=\hat{\tilde{P}}_{\nu} (\hat{\rho}+\mu ),
\end{equation}
where $\hat{\tilde{P}}_\nu (x)= C_+^{\text{\tiny R}}(\nu )
\hat{\tilde{p}}_\nu (x)$ is the dual of
$\tilde{P}_\nu (x)= \hat{C}_+^{\text{\tiny R}}(\nu )
\tilde{p}_\nu (x)$ with the parameters $g_r$
being replaced by the dual parameters $\hat{g}_r$ \eqref{hatgr}.
(For $\nu=0$ the duality relation \eqref{dualpropdeg}
reduces to the evaluation
formula $\tilde{P}_\mu (\rho )=1$ characterizing the
normalization of the polynomials $\tilde{P}_\mu (z )$.)
Just as in the $q$-case, the duality properties for the polynomials
are again inherited by the matrix
$[\tilde{K}_{\mu ,\nu} ]_{\mu ,\nu \in\Lambda_N}$ in \eqref{matrixortho}.
We now have (for parameters
satisfying the truncation condition \eqref{truncdeg})
that transposition of the matrix
$[\tilde{K}_{\mu ,\nu} ]_{\mu ,\nu \in\Lambda_N}$
amounts to the parameter transformation $g_r\rightarrow \hat{g}_r$
and, in particular, that the sum
$\langle 1 , 1\rangle_N^{\text{\tiny R}}$ is invariant with respect to
such a transformation of the parameters.
For parameters satisfying the additional constraint
$g_a-g_b-g_c-g_d=-1$ we have that $\hat{g}_r=g_r$. Hence,
in that case
we are again in a self-dual situation
(i.e. $\hat{\tilde{P}}_\nu (x)= \tilde{P}_\nu (x)$) and the matrix
$[\tilde{K}_{\mu ,\nu} ]_{\mu ,\nu \in\Lambda_N}$ now also becomes
symmetric (in addition to being orthogonal).
\end{remark}

\begin{remark}
It turns out that our multivariable Racah polynomials (i.e., the
multivariable polynomials $\tilde{p}_\lambda (x)$ with parameters
subject to the truncation condition \eqref{truncdeg}) are not the
first generalization of Wilson's one-variable
Racah polynomials to the case of several variables.
Already several years ago Gustafson reported on
a finite system of multivariable orthogonal polynomials
with discrete orthogonality measure tied to
the so-called multiplicity-free Racah coefficients for the group
$U(m+1)$ \cite{gus:whipple's}. In the rank one situation ($m=1$),
the orthogonal
polynomials in question can be reduced to the Racah polynomials
of \cite{wil:some}. It is not clear (at least not to us), though an
interesting question,
whether also for higher rank Gustafson's multivariable Racah polynomials
may be linked with the multivariable Racah polynomials
of the present paper. From the explicit expressions for
the weight functions it seems
that both approaches generalize the one-variable Racah polynomials to
several
variables along very different directions. However, at present we are not
able to rule out completely the possibility
that there might not be some transformation
connecting the two approaches.
Should there indeed exist such a connection (which, however, would
seem more likely after trading the group $U(m+1)$ for $Sp(m)$ say),
then this would imply an interesting link between Gustafson's
group-theoretical program and a degenerate case of the Macdonald theory.
\end{remark}

\begin{remark}\label{trigrem}
The reader may find it illuminating to view the multivariable $q$-Racah
polynomials $p_\lambda$
in Sections~\ref{sec3} and \ref{sec4} as a trigonometric
version of the multivariable Racah polynomials $\tilde{p}_\lambda$
in the present section.
If we substitute Askey-Wilson parameters in accordance with \eqref{repa}
and set $q=e^{i\alpha}$,
then we may rewrite the weight function
$\Delta^{\text{\tiny qR}}$ \eqref{weightqR} as
\begin{equation*}
\Delta^{\text{\tiny qR}} (\nu ) =
\frac{ 1 }
{C_+^{\text{\tiny qR}}(\nu )\: C_-^{\text{\tiny qR}}(\nu ) } ,
\end{equation*}
with
\begin{eqnarray*}
C_+^{\text{\tiny qR}}(\nu )\!\!\! &=&\!\!\! (-4)^{-|\nu |}
\prod_{1\leq j< k\leq n}
\left( \frac{(\rho_j+\rho_k:\sin_\alpha )_{\nu_j+\nu_k}}
            {(g+\rho_j+\rho_k:\sin_\alpha )_{\nu_j+\nu_k}}
\frac{(\rho_j-\rho_k:\sin_\alpha )_{\nu_j-\nu_k}}
     {(g+\rho_j-\rho_k:\sin_\alpha )_{\nu_j-\nu_k}} \right) \\
&& \makebox[2em]{}\times \prod_{1\leq j\leq n}
\left( \frac{(2\rho_j:\sin_\alpha )_{2\nu_j}}
     {\prod_{0\leq r\leq 3}(g_r+\rho_j:\sin_\alpha )_{\nu_j}} \right) ,
\end{eqnarray*}
and
\begin{eqnarray*}
C_-^{\text{\tiny qR}}(\nu )\!\!\! &=&\!\!\! (-4)^{|\nu |} \\
&& \!\!\!\times\!\!\!\prod_{1\leq j< k\leq n}
\left( \frac{(1-g+\rho_j+\rho_k:\sin_\alpha )_{\nu_j+\nu_k}}
            {(1+\rho_j+\rho_k:\sin_\alpha )_{\nu_j+\nu_k}}
       \frac{(1-g+\rho_j-\rho_k:\sin_\alpha )_{\nu_j-\nu_k}}
            {(1+\rho_j-\rho_k:\sin_\alpha )_{\nu_j-\nu_k}} \right) \\
& &\!\!\!\times  \prod_{1\leq j\leq n}
\left( \frac{\prod_{0\leq r\leq 3}(1-g_r+\rho_j:\sin_\alpha )_{\nu_j}}
            {(1+2\rho_j:\sin_\alpha )_{2\nu_j}} \right) \nonumber
\end{eqnarray*}
(where the components of the vector $\rho$ are given by \eqref{rho}).
In the above formula
we have used `trigonometric Pochhammer symbols' defined by
\begin{equation}
(a:\sin_\alpha )_m =
   \sin_\alpha (a) \sin_\alpha (a+1)\cdots \sin_\alpha (a+m-1)
\end{equation}
with $(a:\sin_\alpha )_0\equiv 1$
and $\sin_\alpha (\xi)\equiv \sin (\alpha\xi /2 )$.
Furthermore, we also arrive at corresponding trigonometric expressions for
$\hat{\Delta}^{\text{\tiny qR}}(\lambda)
=1/ (\hat{C}_+^{\text{\tiny qR}}(\lambda)
\hat{C}_-^{\text{\tiny qR}}(\lambda ) )$ and
$\mathcal{N}^{\text{\tiny qR}}(\lambda)=
\hat{C}_-^{\text{\tiny qR}}(\lambda )/
\hat{C}_+^{\text{\tiny qR}}(\lambda)$, which are governed
by dual $c$-functions $\hat{C}_\pm^{\text{\tiny qR}} $
obtained by replacing the parameters $g_r$ by $\hat{g}_r$ \eqref{hatgr}
and the vector $\rho$ by $\hat{\rho}$ \eqref{hatrho}.

It is manifest from these representations that
$\Delta^{\text{\tiny qR}}$,
$\hat{\Delta}^{\text{\tiny qR}}$ and
$\mathcal{N}^{\text{\tiny qR}}$
can be interpreted as trigonometric
versions of
$\Delta^{\text{\tiny R}}$ \eqref{weightR},
$\hat{\Delta}^{\text{\tiny R}}$ (see \eqref{matrixortho}) and
$\mathcal{N}^{\text{\tiny R}}$ \eqref{NR},
respectively.
The transition $q\rightarrow 1$ corresponds to the
limit $\alpha\rightarrow 0$ in which
the renormalized trigonometric Pochhammer symbols
$(2/\alpha)^{m}(a:\sin_\alpha )_m$ go over in the ordinary
Pochhammer symbols $(a)_m$.
As far as the polynomials are concerned, we see that
$p_\lambda (e^{i\alpha x})$ becomes a trigonometric
polynomial in the variables $x_1,\ldots , x_n$
with period $2\pi /\alpha$.
Proposition~\ref{limprp}
describes the rational limit in which the period
of the trigonometric functions tends to infinity (cf. \cite{die:properties})
\begin{equation*}
\tilde{p}_\lambda (x)=
\lim_{\alpha \rightarrow 0}
(i\alpha)^{-2|\lambda |}p_\lambda (e^{i\alpha x}).
\end{equation*}
Notice also that the modified monomial basis elements in \eqref{modbasis}
are in the trigonometric notation of the form
$(-4)^{|\lambda |} \tilde{m}_\lambda (\sin_\alpha (x))$ and
converge, after division by the constants $(i\alpha )^{2|\lambda |}$,
to $\tilde{m}_\lambda (x)$ for $\alpha \rightarrow 0$ .
Finally, in the trigonometric coordinates
the grid points on which the discrete orthogonality measure
$\Delta^{\text{\tiny qR}}$
for the polynomials $p_\lambda (e^{i\alpha x})$
is supported become of the form
$\rho + \nu$, $\nu\in\Lambda_N$.
This means, in particular, that in these coordinates the grid
points do not move when performing the limit $q\rightarrow 1$
(or equivalently $\alpha\rightarrow 0$).
\end{remark}

\section{Positivity domain for $|q|=1$ and the $q$-Racah transform}\label{sec6}
In the preceding sections we have viewed the polynomials and other
objects of interest (such as the discrete weight function
and the difference equation) as rational expressions in the parameters.
As a consequence, we arrived at their properties for generic values
of the parameters. It is clear, however, that the above generic picture
does {\em not} hold for {\em all} values of the parameters.
A most drastic way in which
the generic picture breaks down occurs when the cardinality
of the grid $\tau q^{\nu}$, $\nu\in\Lambda_N$ becomes less than
the dimension of the space $\mathcal{H}^{\text{\tiny qR}}_N$ \eqref{HN}
(or, equivalently, when it becomes
less than the number of points in the alcove $\Lambda_N$
\eqref{alcove}).
This may for instance happen
when $q^M=1$ for $M\in \{ 1,\ldots, N\}$.
The special case $M=1$ (so $q=1$) corresponds of course precisely
to what we already analyzed in more detail in
the previous section by means of a limit transition.
In this section, however, rather than to present any in-depth analysis of
the case that $q$ is such a root of unity
(cf. \cite{che:macdonalds,kir:inner,spi-zhe:zeros}), we will
focus on a parameter domain with $|q|=1$ for which the generic picture
sketched in the preceding sections
{\em does apply} for {\em all} parameter values in the domain and, moreover,
for which the weight function $\Delta^{\text{\tiny qR}}(\nu)$ \eqref{weightqR}
becomes real-valued and positive when $\nu$ lies in the alcove $\Lambda_N$.

To the describe the positivity domain for $|q|=1$
it is convenient to perform the following
trigonometric substitution of the parameters
(cf. Remark~\ref{trigrem})
\begin{equation}\label{trigpa}
\begin{array}{llll}
 q=e^{i\alpha}, &  t=e^{i\alpha g}, & &\\
t_a=e^{i\alpha \mathrm{g}_a}, & t_b = -e^{i\alpha \mathrm{g}_b}, &
t_c =e^{i\alpha (\mathrm{g}_c+1/2)}, & t_d=-e^{i\alpha (\mathrm{g}_d+1/2)} .
\end{array}
\end{equation}
The weights $\Delta^{\text{\tiny qR}}$ \eqref{weightqR} can then
be rewritten with the aid of the trigonometric Pochhammer symbols
(cf. Remark~\ref{trigrem})
\begin{subequations}
\begin{eqnarray}
 (a_1,\ldots a_s:\sin_\alpha )_m
&=&(a_1:\sin_\alpha )_m \cdots (a_s:\sin_\alpha )_m \\
 (a:\sin_\alpha )_m &=&
   \sin_\alpha (a) \sin_\alpha (a+1)\cdots \sin_\alpha (a+m-1) \\
(a_1,\ldots a_s:\cos_\alpha )_m
&=& (a_1:\cos_\alpha )_m \cdots (a_s:\cos_\alpha )_m \\
 (a:\cos_\alpha )_m &=&
   \cos_\alpha (a) \cos_\alpha (a+1)\cdots \cos_\alpha (a+m-1)
\end{eqnarray}
\end{subequations}
where $(a:\sin_\alpha )_0=(a:\cos_\alpha )_0\equiv 1$
and $\sin_\alpha (\xi)\equiv \sin (\alpha\xi /2 )$,
$\cos_\alpha (\xi)\equiv \cos (\alpha\xi /2 )$.
Specifically, we have
\begin{equation}\label{sdweight}
\Delta^{\text{\tiny qR}} (\nu ) =
\frac{ 1 }
{C_+^{\text{\tiny qR}}(\nu )\: C_-^{\text{\tiny qR}}(\nu ) } ,
\end{equation}
with
\begin{eqnarray*}
C_+^{\text{\tiny qR}}(\nu )\!\!\!\!\! &=&\!\!\!\!\!\!\!\!\!
\prod_{1\leq j< k\leq n}\!\!\!
\left( \frac{(\rho_j+\rho_k:\sin_\alpha )_{\nu_j+\nu_k}}
            {(g+\rho_j+\rho_k:\sin_\alpha )_{\nu_j+\nu_k}}
\frac{(\rho_j-\rho_k:\sin_\alpha )_{\nu_j-\nu_k}}
     {(g+\rho_j-\rho_k:\sin_\alpha )_{\nu_j-\nu_k}} \right) \\
&&\hspace{-4.5em}\times \!\!\!\prod_{1\leq j\leq n}\!\!\!
\left( \frac{(\rho_j,1/2+\rho_j:\sin_\alpha )_{\nu_j}
             (\rho_j,1/2+\rho_j:\cos_\alpha )_{\nu_j}  }
     {(\mathrm{g}_a+\rho_j,\mathrm{g}_c+1/2+\rho_j :\sin_\alpha )_{\nu_j}
      (\mathrm{g}_b+\rho_j,\mathrm{g}_d+1/2+\rho_j :\cos_\alpha )_{\nu_j} }
\right) ,
                                                  \nonumber\\
C_-^{\text{\tiny qR}}(\nu )\!\!\!\!\! &=&\!\!\!\!\!\!\!\!\!
\prod_{1\leq j< k\leq n}\!\!\!
\left( \frac{(1-g+\rho_j+\rho_k:\sin_\alpha )_{\nu_j+\nu_k}}
            {(1+\rho_j+\rho_k:\sin_\alpha )_{\nu_j+\nu_k}}
       \frac{(1-g+\rho_j-\rho_k:\sin_\alpha )_{\nu_j-\nu_k}}
            {(1+\rho_j-\rho_k:\sin_\alpha )_{\nu_j-\nu_k}} \right) \\
& &\hspace{-4.5em}\times \!\!\!\prod_{1\leq j\leq n}\!\!\!
\left( \frac{(1-\mathrm{g}_a+\rho_j,1/2-\mathrm{g}_c+\rho_j :\sin_\alpha
)_{\nu_j}
             (1-\mathrm{g}_b+\rho_j,1/2-\mathrm{g}_d+\rho_j :\cos_\alpha
)_{\nu_j}}
            {(1+\rho_j, 1/2+\rho_j:\sin_\alpha )_{\nu_j}
             (1+\rho_j, 1/2+\rho_j:\cos_\alpha )_{\nu_j}} \right) \nonumber
\end{eqnarray*}
and
\begin{equation*}
\rho_j = (n-j) g + \mathrm{g}_a.
\end{equation*}
The corresponding dual objects $\hat{C}_+^{\text{\tiny qR}}$
and $\hat{C}_-^{\text{\tiny qR}}$ are again obtained by
replacing
the parameters $\mathrm{g}_r$ by dual
parameters $\hat{\mathrm{g}}_r$ with
\begin{equation}\label{dualpa}
\left(
\begin{matrix} \hat{\mathrm{g}}_a \\ \hat{\mathrm{g}}_b \\
\hat{\mathrm{g}}_c \\ \hat{\mathrm{g}}_d \end{matrix}
                                                              \right)
= \frac{1}{2} \left(
\begin{array}{rrrr}
1  &  1  &  1  &  1  \\
1  &  1  &  -1  &  -1  \\
1  &  -1  &  1  &  -1  \\
1  &  -1  &  -1  &  1
\end{array}             \right)
\left( \begin{matrix} \mathrm{g}_a \\ \mathrm{g}_b \\
\mathrm{g}_c \\ \mathrm{g}_d \end{matrix} \right)
\end{equation}
and the vector $\rho$ by $\hat{\rho}$ with
$\hat{\rho}_j = (n-j) g + \hat{\mathrm{g}}_a$.
\begin{proposition}\label{posprp}
For parameters satisfying the constraints
\begin{subequations}
\begin{equation}\label{posc1}
\alpha > 0, \;\;\;\;\; g\geq 0 ,\;\;\;\;\;
0\leq \mathrm{g}_a,\mathrm{g}_b < \frac{\pi}{\alpha}, \;\;\;
-\mathrm{g}_a\leq \mathrm{g}_c \leq \mathrm{g}_a,\;\;\;-\mathrm{g}_b\leq
\mathrm{g}_d \leq \mathrm{g}_b
\end{equation}
and the truncation condition
\begin{equation}\label{posc2}
(n-1)g + \mathrm{g}_a +\mathrm{g}_b +N = \frac{\pi}{\alpha} ,
\end{equation}
\end{subequations}
one has that
\begin{eqnarray*}
0< C_\pm^{\text{\tiny qR}} (\nu ) < \infty \;\;\;\;\;\;\; \text{for}\;\;
\nu\in\Lambda_N,\\
0< \hat{C}_\pm^{\text{\tiny qR}} (\mu ) < \infty \;\;\;\;\;\;\; \text{for}\;\;
\mu\in \Lambda_N.
\end{eqnarray*}
\end{proposition}
\begin{proof}
Notice that for $N=0$ the proposition is valid trivially, because
in that case $C_\pm^{\text{\tiny qR}},\hat{C}_\pm^{\text{\tiny qR}}\equiv 1$
(and $\Lambda_N = \{ 0\}$).
Let us from now on assume that $N$ is positive and, furthermore,
let us also temporarily assume that the parameters
$g$ and $\mathrm{g}_a,\mathrm{g}_b$ are nonzero.
It is not very difficult to verify that the conditions \eqref{posc1} and
\eqref{posc2}
then imply that the arguments of the sinus functions in
$C_\pm^{\text{\tiny qR}} (\nu )$ (with $\nu\in \Lambda_N$)
lie in the
open interval $(0, \pi)$ and, similarly,
that the arguments of the cosinus functions
lie in the open interval $(-\pi /2 ,\pi/2)$.
Hence, it follows that the $c$-functions $C_\pm^{\text{\tiny qR}} (\nu )$,
$\nu\in\Lambda_N$ are positive and finite for these parameters.
When one or more of the parameters $g,\mathrm{g}_a,\mathrm{g}_b$ (or
$\mathrm{g}_c,\mathrm{g}_d$
for that matter)
become zero, one may cancel the factors in the numerator/denominator
carrying that parameter
(with the value zero) against the corresponding term
in the denominator/numerator (thus resulting in a trivial unit factor).
By proceeding in this manner one
readily infers that the positivity and finiteness
of the $c$-functions is conserved also
when one or more of the parameters
$g$ and $\mathrm{g}_a,\mathrm{g}_b$ are allowed to become zero.

The dual statement that $\hat{C}_\pm^{\text{\tiny qR}} (\mu )$ is
positive for $\mu\in\Lambda_N$
is now immediate from the observation that the conditions \eqref{posc1},
\eqref{posc2} are self-dual in the sense that they read the same in the
parameters $\mathrm{g}_a,\mathrm{g}_b,\mathrm{g}_c,\mathrm{g}_d$ as in the dual
parameters
$\hat{\mathrm{g}}_a,\hat{\mathrm{g}}_b,\hat{\mathrm{g}}_c,\hat{\mathrm{g}}_d$.
Specifically, the conditions \eqref{posc1}, \eqref{posc2}
(with $N>0$) are equivalent
to the conditions
\begin{subequations}
\begin{equation}\label{posc1d}
\alpha > 0, \;\;\;\;\; g\geq 0 ,\;\;\;\;\;
0\leq \hat{\mathrm{g}}_a,\hat{\mathrm{g}}_b < \frac{\pi}{\alpha}, \;\;\;
-\hat{\mathrm{g}}_a\leq \hat{\mathrm{g}}_c \leq \hat{\mathrm{g}}_a,\;\;\;
-\hat{\mathrm{g}}_b\leq \hat{\mathrm{g}}_d \leq \hat{\mathrm{g}}_b
\end{equation}
and
\begin{equation}\label{posc2d}
(n-1)g + \hat{\mathrm{g}}_a +\hat{\mathrm{g}}_b +N = \frac{\pi}{\alpha} .
\end{equation}
\end{subequations}
That the conditions
\eqref{posc1} and \eqref{posc2} imply
the dual conditions \eqref{posc1d} and \eqref{posc2d}
is seen with the aid of the definition of the
dual parameters in \eqref{dualpa}; that both the conditions
and the dual conditions are
actually equivalent is then clear
from the fact that the parameter transformation determining
the dual parameters in \eqref{dualpa} is an involution.
\end{proof}
\begin{corollary}
For parameters subject to the conditions in Proposition~\ref{posprp},
the weight functions $\Delta^{\text{\tiny qR}}(\nu )$ and
$\hat{\Delta}^{\text{\tiny qR}}(\mu )$
are positive and finite when $\nu,\mu$ lie in the alcove
$\Lambda_N$ \eqref{alcove}.
\end{corollary}
The positivity of Proposition~\ref{posprp} implies
that $p_\lambda$, $\lambda\in\Lambda_N$ is well-defined not
just generically but
for all parameter values in the domain determined
by the conditions \eqref{posc1}, \eqref{posc2}.
\begin{proposition}\label{globalprp}
The polynomials $p_\lambda $, $\lambda\in\Lambda_N$
are well-defined for all parameters \eqref{trigpa}
with values in the domain
determined by the conditions \eqref{posc1} and \eqref{posc2}
(i.e., the expansion coefficients $c_{\lambda ,\mu}$
in \eqref{defcon1} are regular for these parameter values).
\end{proposition}
\begin{proof}
Let us first assume that the parameters are generic (complex say)
but subject to the truncation condition \eqref{trunc}.
For arbitrary $\lambda\in \Lambda_N$, the bilinear form
$\langle \cdot ,\cdot \rangle^{\text{\tiny qR}}_N$ restricted to
the subspace ${\text{Span}} \lbrace m_{\mu} \, | \, \mu<\lambda\rbrace$ is
nondegenerate as a consequence
of Theorem~\ref{orthothm} and Theorem~\ref{sumthm}.
Using the orthogonality of the polynomials and the
nondegeneracy of the bilinear form it is seen that
one can characterize the multivariable
$q$-Racah polynomial corresponding to a
dominant weight vector $\lambda\in\Lambda_N$
as the unique polynomial of the form
\begin{equation*}
p_\lambda(z)= m_{\lambda}(z)+
       \sum_{\mu\in\Lambda_N ,\;\mu<\lambda}c_{\lambda,\mu}m_{\mu}(z)
\end{equation*}
such that
\begin{equation*}
\langle p_\lambda , m_\mu \rangle^{\text{\tiny qR}}_N = 0
\;\;\;\;\;\;\;\;\text{for}\;\; \mu < \lambda .
\end{equation*}
In other words, $p_\lambda$ consists of $m_\lambda$ minus its
(unique) orthogonal
projection with respect to the (nondegenerate) bilinear form
$\langle \cdot ,\cdot \rangle^{\text{\tiny qR}}_N$
onto ${\text{Span}} \lbrace m_{\mu} \, | \, \mu<\lambda\rbrace$.
We thus have the following
inductive Gram-Schmidt-like formula for the (orthogonal)
polynomials $p_\lambda$, $\lambda\in\Lambda_N$
\begin{equation}\label{indfor}
p_\lambda (z) = m_\lambda (z)-
\sum_{\stackrel{\mu\in\Lambda_N}{\mu <\lambda}}
\frac{ \langle m_\lambda ,p_\mu \rangle^{\text{\tiny qR}}_N      }
     {\langle p_\mu , p_\mu \rangle^{\text{\tiny qR}}_N}\; p_\mu (z) .
\end{equation}
(At this point it is helpful again to view Equation~\eqref{indfor}
as a rational identity in the parameters
$q$, $t$ and $t_r$ subject to the truncation condition \eqref{trunc}.)
Using induction on the weight $\lambda$ one sees from this inductive formula
and Theorem~\ref{sumthm} that the polynomial $p_\lambda$ is regular
at parameter values where the $c$-functions
$C_\pm^{\text{\tiny qR}}, \hat{C}_\pm^{\text{\tiny qR}}$ are
regular and nonzero.
The proposition is then immediate from Proposition~\ref{posprp}.
\end{proof}
With the aid of Proposition~\ref{posprp} and \ref{globalprp}
it is seen that the results of Section~\ref{sec3} and \ref{sec4}
hold for all parameter values in the domain determined
by the conditions \eqref{posc1} and \eqref{posc2}. At some points the
positivity of the measure will enable us to even formulate
a somewhat stronger version of the results stated there.
For instance, for the parameters in the positivity domain \eqref{posc1},
\eqref{posc2} the polynomial $p_\lambda (e^{i\alpha x})$  ($x\in\mathbb{R}^n$)
and thus also
the matrix $[K_{\mu ,\nu}]_{\mu ,\nu\in\Lambda_N}$ \eqref{qmatrixortho}
is real. Hence, in addition to being orthogonal the matrix
$[K_{\mu ,\nu}]_{\mu ,\nu\in\Lambda_N}$ now also becomes unitary.
(The real-valuedness of polynomials $p_\lambda (e^{i\alpha x})$
for parameters in the positivity domain \eqref{posc1}, \eqref{posc2}
follows e.g. from the inductive formula \eqref{indfor}
together with the observation that the (even) monomials
$m_\lambda (e^{i\alpha x})$ are real.)
The fact that both polynomials and weight function
are real for parameters in the positivity domain \eqref{posc1}, \eqref{posc2}
allows us to restrict the bilinear form
$\langle \cdot ,\cdot \rangle^{\text{\tiny qR}}_N$ \eqref{discreteIP}
to a real form (by restricting to the real vector space spanned
by $m_\lambda$, $\lambda\in\Lambda_N$), which in turn
can be extended to a positive definite sesquilinear
form on the complex vector space $\mathcal{H}_N^{\text{\tiny qR}}$ \eqref{HN}
(in the standard way).

In the remainder of this section we will interpret the orthogonality
and orthonormalization properties of the multivariable $q$-Racah
polynomials with parameters in the positivity domain
in terms of a finite-dimensional discrete integral transformation
for grid functions.
To this end we need to introduce some further notation.
Let $L^2(\rho +\Lambda_N, \Delta^{\text{\tiny qR}} )$ be
the finite-dimensional Hilbert space of complex functions over
the grid points $\rho +\nu$, $\nu\in\Lambda_N$
endowed with the standard inner product
determined by the positive weights $\Delta^{\text{\tiny qR}}(\nu)$,
$\nu\in\Lambda_N$
\begin{equation}\label{sesqui}
\langle f , g\rangle_{\Delta } = \sum_{\nu \in\Lambda_N}
f(\rho +\nu )\, \overline {g(\rho +\nu )}\:
\Delta^{\text{\tiny qR}}(\nu).
\end{equation}
Similarly, the space
$L^2(\hat{\rho} +\Lambda_N, \hat{\Delta}^{\text{\tiny qR}} )$ denotes the
corresponding dual Hilbert space consisting of the complex functions over
the grid points $\hat{\rho} +\mu$, $\mu\in\Lambda_N$
endowed with the standard inner product
$\langle \cdot , \cdot\rangle_{\hat{\Delta} }$
determined by the positive weights $\hat{\Delta}^{\text{\tiny qR}}(\mu)$,
$\mu\in\Lambda_N$.
We define the operator
$\mathcal{K}:L^2(\rho +\Lambda_N, \Delta^{\text{\tiny qR}} )
\longrightarrow L^2(\hat{\rho} +\Lambda_N, \hat{\Delta}^{\text{\tiny qR}} )$
as the map with kernel
\begin{equation}\label{kern}
\mathcal{K}(\hat{\rho}+\mu ,\rho +\nu) \equiv
\frac{P_\mu (e^{i\alpha (\rho +\nu)} ) \Delta^{\text{\tiny qR}} (\nu ) }
     {\sqrt{\langle 1 ,1\rangle_N^{\text{\tiny qR}}}}
\end{equation}
i.e., the operator $\mathcal{K}$ acts on a grid function
$f: \rho +\Lambda_N\rightarrow \mathbb{C}$ as
\begin{equation}
(\mathcal{K}f)(\hat{\rho}+\mu)=
\sum_{\nu\in\Lambda_N}\mathcal{K}(\hat{\rho}+\mu ,\rho +\nu)f(\rho +\nu )
\end{equation}
(thus producing a function
$(\mathcal{K}f):\hat{\rho} +\Lambda_N\rightarrow \mathbb{C}$).
Here we have employed the renormalized multivariable Askey-Wilson/$q$-Racah
polynomials $P_\lambda (z)= \hat{C}^{\text{\tiny qR}}_+(\lambda) p_\lambda (z)$
satisfying the normalization condition
$P_\lambda (\tau )=1$ (cf. Remark~\ref{dualrem}).
We also define the dual map
$\hat{\mathcal{K}}:L^2(\hat{\rho} +\Lambda_N, \hat{\Delta}^{\text{\tiny qR}} )
\longrightarrow L^2(\rho +\Lambda_N, \Delta^{\text{\tiny qR}} )$ determined by
the kernel
\begin{equation}\label{dkern}
\hat{\mathcal{K}}(\rho+\nu ,\hat{\rho} +\mu) \equiv
\frac{\hat{P}_\nu (e^{i\alpha (\hat{\rho} +\mu)} )
      \hat{\Delta}^{\text{\tiny qR}} (\mu )}
     {\sqrt{\langle 1 ,1\rangle_N^{\text{\tiny qR}}}}
\end{equation}
($\hat{P}_\lambda (z) = C^{\text{\tiny qR}}_+(\lambda) \hat{p}_\lambda (z)$)
and acting on a
function $\hat{f}: \hat{\rho} +\Lambda_N\rightarrow \mathbb{C}$
over the dual grid as
\begin{equation}
(\hat{\mathcal{K}}\hat{f})(\rho+\nu)=
\sum_{\mu\in\Lambda_N}\hat{\mathcal{K}}(\rho+\nu ,\hat{\rho} +\mu)
\hat{f}(\hat{\rho} +\mu ).
\end{equation}

The following theorem describes a discrete integral
transform---the `$q$-Racah transform'---between
grid functions in
$L^2(\rho +\Lambda_N, \Delta^{\text{\tiny qR}} )$ and
$L^2(\hat{\rho} +\Lambda_N, \hat{\Delta}^{\text{\tiny qR}} )$
together with its inversion formula.
\begin{theorem}\label{trafothm}
For parameters subject to the conditions \eqref{posc1}, \eqref{posc2}
the map
$\mathcal{K}:L^2(\rho +\Lambda_N, \Delta^{\text{\tiny qR}} )
\longrightarrow L^2(\hat{\rho} +\Lambda_N, \hat{\Delta}^{\text{\tiny qR}} )$
is an isometric isomorphism. The inverse of ${\mathcal{K}}$ is given by the map
$\hat{\mathcal{K}}:L^2(\hat{\rho} +\Lambda_N, \hat{\Delta}^{\text{\tiny qR}} )
\longrightarrow L^2(\rho +\Lambda_N, \Delta^{\text{\tiny qR}} )$.
\end{theorem}
\begin{proof}
Let $\mathcal{K}=[\mathcal{K}_{\mu ,\nu}]_{\mu ,\nu \in\Lambda_N}$
be the matrix with elements
$\mathcal{K}_{\mu ,\nu}\equiv \mathcal{K}(\hat{\rho}+\mu ,\rho +\nu )=
P_\mu (e^{i\alpha (\rho +\nu)} ) \Delta^{\text{\tiny qR}} (\nu ) /
(\langle 1 ,1\rangle_N^{\text{\tiny qR}})^{1/2}$ and
let $\Delta$ and $\hat{\Delta}$ be the diagonal matrices with the quantities
$\Delta^{\text{\tiny qR}}(\nu )$ and $\hat{\Delta}^{\text{\tiny qR}}(\mu )$
($\nu ,\mu \in \Lambda_N$) on the diagonal, respectively.
(Here it is of course
again assumed that the columns and rows are ordered
by a total extension of the partial order \eqref{order}.)
We then have that
\begin{equation*}
\mathcal{K}=\hat{\Delta}^{-1/2} K \Delta^{1/2},
\end{equation*}
where $K=[K_{\mu ,\nu}]_{\mu ,\nu \in\Lambda_N}$ denotes the matrix
given in \eqref{qmatrixortho}.
The unitarity of the (real orthogonal) matrix
$K$ (for parameters in the positivity domain)
and the inversion formula $K^{-1}=\hat{K}$ imply for
the matrix $\mathcal{K}$ that
\begin{equation}\label{matKprop}
\mathcal{K}^{*}\hat{\Delta}\mathcal{K}=\Delta
\;\;\;\;\;\;\text{and that}\;\;\;\;\;\;
\mathcal{K}^{-1}=\hat{\mathcal{K}}
\end{equation}
(where $\hat{\mathcal{K}}$ and
$\mathcal{K}^{*}$ are the dual and the adjoint (transpose) of $\mathcal{K}$).
Clearly, Formula~\eqref{matKprop}
boils down to a reformulation of the statement in the theorem.
(The first property in \eqref{matKprop} says that the matrix
$\mathcal{K}$ determines an isometry
between the inner product spaces
endowed with the positive sesquilinear
forms associated to $\Delta$ and $\hat{\Delta}$, respectively.)
\end{proof}
It is instructive to view the $q$-Racah transform as a Fourier type
transformation between the grid functions in
$L^2(\rho +\Lambda_N, \Delta^{\text{\tiny qR}} )$ and
$L^2(\hat{\rho} +\Lambda_N, \hat{\Delta}^{\text{\tiny qR}} )$:
\begin{subequations}
\begin{eqnarray}
\overset{\wedge}{f}(\hat{\rho}+\mu )\equiv\;\;\;
(\mathcal{K}f)(\hat{\rho}+\mu )&=&
\frac{\langle f ,P_\mu\rangle_\Delta}{\sqrt{\langle 1,1\rangle_\Delta} } , \\
f(\rho+\nu )= (\mathcal{K}^{-1}\overset{\wedge}{f})(\rho+\nu )&=&
\frac{\langle \overset{\wedge}{f} ,\hat{P}_\nu \rangle_{\hat{\Delta}}}
     {\sqrt{\langle 1,1\rangle_{\hat{\Delta}} }} ,
\end{eqnarray}
\end{subequations}
where $\langle \cdot ,\cdot \rangle_\Delta$ is taken from
\eqref{sesqui} and
$\langle \cdot ,\cdot \rangle_{\hat{\Delta}}$ denotes its dual version
(and the functions $P_\mu$ and $\hat{P}_\nu$ stand for
the restrictions to the grids of
$P_\mu (e^{i\alpha x})$ and $\hat{P}_\nu (e^{i\alpha x})$, respectively).
(Recall also that for parameters in the positivity domain
\eqref{posc1}, \eqref{posc2} one has that
$\langle 1,1\rangle_{\hat{\Delta}} =\langle 1,1\rangle_{\Delta} $
as a consequence of Remark~\ref{dualrem}.)

We will next discuss the behavior of a discretization of the difference
operator $D$ \eqref{DAW} with respect to the $q$-Racah transform.
Let $D^{\text{\tiny qR}}:
L^2(\rho +\Lambda_N, \Delta^{\text{\tiny qR}} )\longrightarrow
L^2(\rho +\Lambda_N, \Delta^{\text{\tiny qR}} )$ be
the discrete difference operator of the form
\begin{equation}\label{Ddis}
D^{\text{\tiny qR}} =
\sum_{\stackrel{1\leq j\leq n}
               {\nu +e_j\in\Lambda_N }}
V_{+j}(\rho+ \nu) (T_j -1) +
\sum_{\stackrel{1\leq j\leq n}
               {\nu -e_j\in\Lambda_N }}
V_{-j}(\rho+ \nu) (T_j^{-1} -1)
\end{equation}
where
\begin{equation*}
V_{\pm j}(x)= w(\pm x_j)
\prod_{\stackrel{1\leq k\leq n}{k\neq j}} v(\pm x_j+x_k)\, v(\pm x_j-x_k),
\end{equation*}
with
\begin{eqnarray*}
v(\xi ) &=& \frac{\sin \frac{\alpha }{2}(g+\xi )}
                   {\sin (\frac{\alpha \xi }{2})}, \\
w(\xi ) &=&  \frac{\sin \frac{\alpha}{2}(\mathrm{g}_a+\xi )}
                   {\sin (\frac{\alpha \xi }{2})}
               \frac{\cos \frac{\alpha}{2}(\mathrm{g}_b+\xi )}
                   {\cos (\frac{\alpha \xi }{2})}
               \frac{\sin \frac{\alpha}{2}(\mathrm{g}_c+1/2+\xi )}
                   {\sin \frac{\alpha}{2} (1/2+\xi )}
               \frac{\cos \frac{\alpha}{2}(\mathrm{g}_d+1/2+\xi )}
                   {\cos \frac{\alpha}{2} (1/2+ \xi )}
\end{eqnarray*}
and the action of the operators $T_j^{\pm 1}$ is given by
\begin{equation*}
(T_j^{\pm 1} f) (\rho +\nu ) = f(\rho + \nu \pm e_j).
\end{equation*}
Notice that the conditions $\nu +e_j \in\Lambda_N$ and $\nu -e_j\in \Lambda_N$
in the summations of \eqref{Ddis}
guarantee that the function $(D^{\text{\tiny qR}}f)(\rho +\nu )$ for
$\nu\in\Lambda_N$ indeed depends only on the values of $f$ on the grid points
in $\rho +\Lambda_N$.
Up to (multiplication by) an overall constant factor with value
$t^{n-1}(t_0t_1t_2t_3q^{-1})^{1/2}=
e^{i\alpha [ (n-1)g +(\mathrm{g}_0+\cdots +\mathrm{g}_3)/2]}$
the operator
$D^{\text{\tiny qR}}$ \eqref{Ddis} amounts to the restriction of
the operator $D$ \eqref{DAW}
to grid functions (cf. Remark~\ref{restrem}) rewritten in trigonometric
form (recall \eqref{trigpa}).
Such a restriction is well-defined because the coefficients of
the difference operator $D$ \eqref{DAW} are regular on the grid points
and they vanish (cf. Lemma~\ref{reslem})
when the $q$-shift operators $T_{j,q}^{\pm 1}$
in \eqref{DAW} shift
the argument of a function out of the grid
(which leads us to the above-mentioned
restrictions on the sums in the discretized
operator \eqref{Ddis}).
Let us
furthermore introduce the multiplication operator
$E:
L^2(\rho +\Lambda_N,\Delta^{\text{\tiny qR}} )\rightarrow
L^2(\rho +\Lambda_N,\Delta^{\text{\tiny qR}} )$
defined by
\begin{equation}
(E f) (\rho +\nu ) = E_\nu  \; f(\rho +\nu )
\end{equation}
with
\begin{equation*}
E_\nu = 2\sum_{1\leq j\leq n} \biggl( \cos \alpha (\rho_j +\nu_j )-
                     \cos (\alpha \rho_j) \biggr),
\end{equation*}
together with the corresponding dual operators
$\hat{D}^{\text{\tiny qR}}:
L^2(\hat{\rho} +\Lambda_N, \hat{\Delta}^{\text{\tiny qR}} )\longrightarrow
L^2(\hat{\rho} +\Lambda_N, \hat{\Delta}^{\text{\tiny qR}} )$
and $\hat{E}:
L^2(\hat{\rho} +\Lambda_N,\hat{\Delta}^{\text{\tiny qR}} )\rightarrow
L^2(\hat{\rho} +\Lambda_N,\hat{\Delta}^{\text{\tiny qR}} )$
in which $\mathrm{g}_r\rightarrow \hat{\mathrm{g}}_r$ and $\rho\rightarrow
\hat{\rho}$.
(The dual quantities $\hat{E}_\lambda= 2\sum_{1\leq j\leq n}
( \cos \alpha (\hat{\rho}_j +\lambda_j )- \cos (\alpha \hat{\rho}_j) )$
coincide with the eigenvalues $E_{\lambda ,\lambda }$ in \eqref{eigenvAW}
up to multiplication by the constant factor
$t^{n-1}(t_0t_1t_2t_3q^{-1})^{1/2}=
e^{i\alpha [ (n-1)g +(\mathrm{g}_0+\cdots +\mathrm{g}_3)/2]}$
relating $D^{\text{\tiny qR}}$ \eqref{Ddis} and $D$ \eqref{DAW} .)

The following theorem states that the $q$-Racah transform (i.e the
discrete integral transformation $\mathcal{K}$)
is the eigenfunction transformation that diagonalizes the
operators $D^{\text{\tiny qR}}$ and $\hat{D}^{\text{\tiny qR}}$.
\begin{theorem}\label{diagonalthm}
For parameters in the positivity domain \eqref{posc1},
\eqref{posc2} the discrete difference operators
$D^{\text{\tiny qR}}:L^2(\rho +\Lambda_N, \Delta^{\text{\tiny qR}} )
\longrightarrow
L^2(\rho +\Lambda_N, \Delta^{\text{\tiny qR}} )$
and $\hat{D}^{\text{\tiny qR}}:
L^2(\hat{\rho} +\Lambda_N,\hat{\Delta}^{\text{\tiny qR}} )\longrightarrow
L^2(\hat{\rho} +\Lambda_N,\hat{\Delta}^{\text{\tiny qR}} )$
are self-adjoint and
the map $\mathcal{K}:L^2(\rho +\Lambda_N, \Delta^{\text{\tiny qR}} )
\longrightarrow L^2(\hat{\rho} +\Lambda_N, \hat{\Delta}^{\text{\tiny qR}} )$
constitutes the corresponding (unitary)
eigenfunction transformation diagonalizing these
operators
\begin{equation}\label{diagonalization}
\mathcal{K} D^{\text{\tiny qR}} \mathcal{K}^{-1} = \hat{E},
\;\;\;\;\;\;\;\;\;\;\;\;\;\;\;\;\;\;\;
\mathcal{K}^{-1} \hat{D}^{\text{\tiny qR}} \mathcal{K}= E.
\end{equation}
\end{theorem}
\begin{proof}
Clearly it is sufficient to prove only one of the diagonalization formulas
in \eqref{diagonalization} because the other will then automatically follow
upon dualization.
Since for parameters in the positivity domain
the (real) functions
$P_\lambda (e^{i\alpha (\rho +\nu)})$,
$\lambda\in\Lambda_N$ form
an orthogonal basis for the space
$L^2(\rho +\Lambda_N, \Delta^{\text{\tiny qR}} )$, it is enough
to show that
\begin{equation}\label{onbase}
(\mathcal{K} D^{\text{\tiny qR}}) P_\lambda =
(\hat{E}\mathcal{K}) P_\lambda ,\;\;\;\;\;\;\;\text{for all}\;\;\;
\lambda\in\Lambda_N
\end{equation}
(where $P_\lambda$ stands for the discretized
trigonometric polynomial $P_\lambda (e^{i\alpha (\rho +\nu)})$).
Equation~\eqref{onbase} is immediate from the discretized eigenvalue equation
(cf. Remark~\ref{restrem})
\begin{equation*}
D^{\text{\tiny qR}}P_\lambda =\hat{E}_\lambda P_\lambda
\end{equation*}
and the observation that $(\mathcal{K}P_\lambda )(\hat{\rho} +\mu )$ is
nonzero only in the point $\hat{\rho} +\lambda$ (i.e. for $\mu =\lambda$)
in view of the orthogonality of the multivariable $q$-Racah polynomials.
(We thus have that $\mathcal{K} D^{\text{\tiny qR}} P_\lambda=
\hat{E}_\lambda \mathcal{K} P_\lambda =\hat{E}\mathcal{K} P_\lambda$.)
The self-adjointness of $D^{\text{\tiny qR}}$ and
$\hat{D}^{\text{\tiny qR}}$ now follows from
the fact that the discrete difference operators
are unitarily equivalent to the real (for parameters
in the positivity domain) multiplication operators $\hat{E}$
and $E$, respectively.
\end{proof}

\begin{remark}
If the parameters satisfy the additional constraint
\begin{equation}
\mathrm{g}_a-\mathrm{g}_b-\mathrm{g}_c-\mathrm{g}_d=0,
\end{equation}
then $\hat{\mathrm{g}}_r =\mathrm{g}_r$, $r=0,\ldots ,3$ (see \eqref{dualpa}).
Hence, we are then in a self-dual situation (cf. Remark~\ref{dualrem})
with both the Hilbert spaces
$L^2(\rho +\Lambda_N, \Delta^{\text{\tiny qR}} )$ and
$L^2(\hat{\rho} +\Lambda_N, \hat{\Delta}^{\text{\tiny qR}} )$ coinciding
and the map $\mathcal{K}:L^2(\rho +\Lambda_N, \Delta^{\text{\tiny qR}} )
\longrightarrow L^2(\rho +\Lambda_N, \Delta^{\text{\tiny qR}} )$
being an involution ($\mathcal{K}^2 = \hat{\mathcal{K}}\mathcal{K}=Id$).
\end{remark}

\begin{remark}
A positivity domain for a self-dual
one-parameter subfamily of the
one-variable $q$-Racah polynomials with $|q|=1$ similar to the domain
considered in this section
can be found in Section~3C2 of \cite{rui:finite} together
with a discussion of the corresponding
finite-dimensional discrete integral transform.
\end{remark}

\begin{remark}
The trigonometric polynomials
$p_\lambda (e^{i\alpha x})$ are
invariant with respect to permutations, sign flips ($x_j\rightarrow -x_j$)
and translations ($x_j\rightarrow x_j+2\pi/\alpha$)
of the variables $x_1,\ldots ,x_n$.
A fundamental domain for $\mathbb{R}^n$ modulo the action of the
discrete symmetry group generated by the permutations, sign flips
and translations over a period $2\pi /\alpha$
is given by the (Weyl) alcove
\begin{equation}\label{Walcove}
\{ x\in \mathbb{R}^n | \; \pi/\alpha \geq x_1\geq x_2
   \geq \cdots \geq x_n\geq 0 \} .
\end{equation}
It is therefore natural to consider the trigonometric polynomials
$p_\lambda (e^{i\alpha x})$, $\lambda\in\Lambda$ as polynomials
over the alcove \eqref{Walcove}.
In the trigonometric context the necessity to truncate the grid
with points $\rho +\nu$, $\nu\in\Lambda$ on which the
masses of the discrete orthogonality measures are concentrated arises
as a natural consequence of the periodicity of the trigonometric
functions (which demands that the support of
the discrete orthogonality measure be finite).
The conditions \eqref{posc1}, \eqref{posc2} in Proposition~\ref{posprp}
arrange things in such a manner
that $\Delta^{\text{\tiny qR}}(\nu )$
\eqref{sdweight} is positive for $\nu\in\Lambda_N$ and zero
for $\nu\in\Lambda\setminus\Lambda_N$, and guarantee furthermore
that the grid $\rho +\Lambda_N$ supporting
the discrete orthogonality measure for the polynomials
$p_\lambda (e^{i\alpha x})$, $\lambda\in\Lambda_N$ fits
in the fundamental domain \eqref{Walcove}.
\end{remark}

\begin{remark}
The second-order operator $D$ \eqref{DAW} sits in a commutative algebra
generated by $n$-independent commuting analytic
difference operators $D_1,D_2,\ldots ,D_n$ of order $2,4,\ldots ,2n$,
respectively \cite{die:commuting,die:self-dual}.
The first of these operators, viz. $D_1$, corresponds to the
Koornwinder-Macdonald difference
operator $D$ of Section~\ref{sec2}.
After restriction to functions on the grid $\rho+\Lambda_N$
the analytic difference operators
go over in a family of (commuting) discrete difference operators
$D_r^{\text{\tiny qR}}:L^2(\rho +\Lambda_N, \Delta^{\text{\tiny qR}})
\longrightarrow L^2(\rho +\Lambda_N, \Delta^{\text{\tiny qR}})$
($r=1,\ldots ,n$) given explicitly by (cf. Remark~\ref{restrem})
\begin{equation}\label{ddos}
D_r^{\text{\tiny qR}}=
\sum_{\stackrel{J\subset \{ 1,\ldots ,n\} ,\, 0\leq|J|\leq r}
               {\varepsilon_j=\pm 1,\; j\in J;\;
                e_{\varepsilon J} +\nu \in \Lambda_N}}
   \!\!\!\!
U_{J^c,\, r-|J|}(\rho +\nu )\,  V_{\varepsilon J,\, J^c}(\rho +\nu )\,
T_{\varepsilon J },\;\;\;\;\;\;\;\;\;
r=1,\ldots ,n,
\end{equation}
with $(T_{\varepsilon J }f)(\rho +\nu )=f(\rho +\nu +e_{\varepsilon J })$,
$e_{\varepsilon J }=\sum_{j\in J} \varepsilon_j e_j$, and
\begin{eqnarray*}
V_{\varepsilon J,\, K}(x) &=&
\prod_{j\in J} w(\varepsilon_jx_j)
\prod_{\stackrel{ j,j^\prime \in J}{j<j^\prime}}
v(\varepsilon_jx_j+\varepsilon_{j^\prime}x_{j^\prime})
 v(\varepsilon_jx_j+\varepsilon_{j^\prime}x_{j^\prime}+1 )\\
& & \times  \prod_{\stackrel{j\in J}{k\in K}} v(\varepsilon_j x_j+x_k)
v(\varepsilon_j x_j -x_k),  \\ [1ex]
U_{K,p}(x)&=& (-1)^p \!\! \sum_{\stackrel{L\subset K,\, |L|=p}
                               {\varepsilon_l =\pm 1,\; l\in L}} \;
\prod_{l\in L} w(\varepsilon_l x_l)\,
\prod_{\stackrel{l,l^\prime \in L}{l<l^\prime}}
v(\varepsilon_lx_l+\varepsilon_{l^\prime}x_{l^\prime})
v(-\varepsilon_lx_l-\varepsilon_{l^\prime}x_{l^\prime}-1 ) \\
& & \times \prod_{\stackrel{l\in L}{k\in K\setminus L}}
v(\varepsilon_l x_l+x_k) v(\varepsilon_l x_l -x_k)
\end{eqnarray*}
(where $v$ and $w$ are the same as in \eqref{Ddis}).
In the above expressions
we have used the conventions that empty products are equal to one, and
that $U_{K,p}\equiv 1$ for $p=0$.
Notice that the coefficient functions
$V_{\varepsilon J,\, K}(x)$ and $U_{K,p}(x)$
are regular for $x\in \rho +\Lambda_N$ and that
the condition $e_{\varepsilon J}+\nu \in \Lambda_N$
in the summation
again guarantees that $(D_r^{\text{\tiny qR}}f)(\rho +\nu )$
only depends on the values of $f$ on the grid points in $\rho +\Lambda_N$
when $\nu$ lies in $\Lambda_N$. (Hence the operator
$D_r^{\text{\tiny qR}}$ \eqref{ddos} is well-defined
as an operator in $L^2(\rho +\Lambda_N, \Delta^{\text{\tiny qR}})$.)
After introducing also the corresponding multiplication operators
$E_r:L^2(\rho +\Lambda_N, \Delta^{\text{\tiny qR}})\longrightarrow
L^2(\rho +\Lambda_N, \Delta^{\text{\tiny qR}})$ given by
\begin{equation*}
(E_r f)(\rho +\nu ) = E_{r, \nu } f(\rho +\nu )
\end{equation*}
($r=1,\ldots ,n$) with
\begin{eqnarray*}
E_{r,\nu}&=& 2^r\!\!\!\sum_{\stackrel{J\subset \{
 1,\ldots ,n\} }{0\leq |J|\leq r}} \!\! (-1)^{r-|J|} \Bigl(
 \prod_{j\in J} \cos\alpha (\rho_j+\nu_j) \\
& & \makebox[7em]{} \times \sum_{r\leq
 l_1\leq\cdots\leq l_{r-|J|}\leq n}\!\!\!\!  \cos (\alpha
 \rho_{l_1})\cdots \cos (\alpha \rho_{l_{r-|J|}}) \Bigr)
\end{eqnarray*}
(where the second sum in $E_{r,\nu}$ should be read as $1$ when $|J|=r$)
and the associated dual difference operators
$\hat{D}_r^{\text{\tiny qR}}:L^2(\hat{\rho} +\Lambda_N,
\hat{\Delta}^{\text{\tiny qR}})\longrightarrow
L^2(\hat{\rho} +\Lambda_N, \hat{\Delta}^{\text{\tiny qR}})$ and
dual multiplication operators $\hat{E}_{r}:L^2(\hat{\rho} +\Lambda_N,
\hat{\Delta}^{\text{\tiny qR}})\longrightarrow
L^2(\hat{\rho} +\Lambda_N, \hat{\Delta}^{\text{\tiny qR}})$,
we are in the position to formulate a
generalization of Theorem~\ref{diagonalthm} pertaining to
these higher-order discrete difference operators.
\begin{theorem}\label{hodisDthm}
For parameters in the positivity domain \eqref{posc1},
\eqref{posc2} the commuting discrete difference operators
$D_r^{\text{\tiny qR}}:L^2(\rho +\Lambda_N, \Delta^{\text{\tiny qR}} )
\longrightarrow
L^2(\rho +\Lambda_N, \Delta^{\text{\tiny qR}} )$
and $\hat{D}_r^{\text{\tiny qR}}:
L^2(\hat{\rho} +\Lambda_N,\hat{\Delta}^{\text{\tiny qR}} )\longrightarrow
L^2(\hat{\rho} +\Lambda_N,\hat{\Delta}^{\text{\tiny qR}} )$
are self-adjoint and
the map $\mathcal{K}:L^2(\rho +\Lambda_N, \Delta^{\text{\tiny qR}} )
\longrightarrow L^2(\hat{\rho} +\Lambda_N, \hat{\Delta}^{\text{\tiny qR}} )$
constitutes the (unitary) joint
eigenfunction transformation simultaneously
diagonalizing these operators
\begin{equation*}
\mathcal{K} D_r^{\text{\tiny qR}} \mathcal{K}^{-1} = \hat{E}_r,
\;\;\;\;\;\;\;\;\;\;\;\;\;\;\;\;\;\;\;
\mathcal{K}^{-1} \hat{D}_r^{\text{\tiny qR}} \mathcal{K}= E_r
\end{equation*}
($r=1,\ldots ,n$).
\end{theorem}
The proof of Theorem~\ref{hodisDthm}
runs along the same lines as that of Theorem~\ref{diagonalthm} and
hinges on the discrete difference equations in Remark~\ref{restrem} and the
real-valuedness of the multiplication operators $E_r$ and $\hat{E}_r$
for parameters in the positivity domain \eqref{posc1}, \eqref{posc2}.
For $r=1$ Theorem~\ref{hodisDthm} reduces to Theorem~\ref{diagonalthm}.
\end{remark}

\appendix
\section{Proof for the symmetry of $D$}\label{appa}
In this appendix we prove Proposition~\ref{symprp}, which stated
that the $q$-difference operator $D$ \eqref{DAW} is symmetric
with respect to the bilinear form
$\langle \cdot ,\cdot \rangle_N^{\text{\tiny qR}}$ \eqref{discreteIP},
for parameters
satisfying the truncation condition \eqref{trunc}.
This proposition was a key ingredient in our orthogonality proof for the
multivariable $q$-Racah polynomials (i.e. the multivariable Askey-Wilson
polynomials with parameters subject to the truncation condition)
with respect to the bilinear form
$\langle \cdot ,\cdot \rangle_N^{\text{\tiny qR}}$ \eqref{discreteIP}
(Theorem~\ref{orthothm}).

Let us first recall the explicit form of the
functions $V_{\pm j}(z)$ that determine the coefficients of the
difference operator $D$ \eqref{DAW}:
\begin{eqnarray}\label{coefV}
V_{\varepsilon j}(z)&=&
\frac{(1-t_0 z_j^{\varepsilon})
      (1-t_1 z_j^{\varepsilon})
      (1-t_2 z_j^{\varepsilon})
      (1-t_3 z_j^{\varepsilon})}
     {(1-z_j^{2\varepsilon})(1-q z_j^{2\varepsilon})} \\
& &\makebox[2em]{}\times \prod_{\stackrel{1\leq k\leq n}{k\neq j}}
\frac{(1-t\, z_j^{\varepsilon}\, z_k)(1-t\, z_j^{\varepsilon}\, z_k^{-1})}
     {(1-z_j^{\varepsilon}\, z_k)(1-z_j^{\varepsilon}\,z_k^{-1}) },
   \;\;\;\;\;\;\; \varepsilon =\pm 1. \nonumber
\end{eqnarray}
The symmetry proof for the operator $D$ hinges on two lemmas.
The first lemma describes a relation between the
discrete weight function
$\Delta^{\text{\tiny qR}}(\nu)$ \eqref{weightqR}
and the coefficients
$V_{\varepsilon j}(z)$ of $D$ evaluated
at the grid points $\tau q^\nu$ \eqref{grid}.

\begin{lemma}\label{fliplem}
Let us assume that $\nu$ and $\nu +\varepsilon e_j$ are in the dominant cone
$\Lambda$ \eqref{cone} (here $\varepsilon$ is $+1$ or $-1$).
Then
\begin{equation}\label{relVweight}
\Delta^{\text{\tiny qR}}(\nu +\varepsilon e_j)
V_{-\varepsilon j}(\tau q^{\nu +\varepsilon e_j}) =
\Delta^{\text{\tiny qR}} (\nu ) V_{\varepsilon j}(\tau q^\nu )
\end{equation}
(where $\tau q^\nu$, $\Delta^{\text{\tiny qR}}(\nu)$ and
$V_{\varepsilon j}(z)$
are given by \eqref{grid}, \eqref{weightqR} and \eqref{coefV}).
\end{lemma}
\begin{proof}
The relation \eqref{relVweight} between
$\Delta^{\text{\tiny qR}}=
1/(C^{\text{\tiny qR}}_+\, C^{\text{\tiny qR}}_-)$
\eqref{weightqR} and $V_{\pm j}$
follows from the difference equations
\begin{equation}\label{CV+}
 \frac{C^{\text{\tiny qR}}_+(\tilde{\nu} -e_j)}
     {C^{\text{\tiny qR}}_+(\tilde{\nu})}
    = V_{+j}(\tau q^{\tilde{\nu} -e_j})f_j(\tilde{\nu})
\;\;\;\;\;\;\;\;\;\;\;\;\;\;\;\;\;\;\;\;
\text{for}\;\;\;\;\tilde{\nu}, \tilde{\nu}-e_j \in \Lambda
\end{equation}
and
\begin{equation}\label{CV-}
\frac{C^{\text{\tiny qR}}_-(\tilde{\nu} +e_j)}
        {C^{\text{\tiny qR}}_-(\tilde{\nu})}
    = V_{-j}(\tau q^{\tilde{\nu} +e_j})f_j(\tilde{\nu}+e_j)
\;\;\;\;\;\;\;\;\;\;\;\;\;
\text{for}\;\;\;\; \tilde{\nu}, \tilde{\nu}+e_j \in \Lambda
\end{equation}
linking the $c$-functions $C^{\text{\tiny qR}}_\pm$
to the coefficients $V_{\pm j}$.
The factors $f_j$ in \eqref{CV+} and \eqref{CV-}
represent certain intermediate products of the form
\begin{equation*}
f_j(\tilde{\nu})\equiv t^{j-n}(t_0t_1t_2t_3q^{-1})^{-1/2}
\prod_{1\leq k <j}
\frac{(1-\tau_j^{-1}\tau_kq^{\tilde{\nu}_k-\tilde{\nu}_j})
(1-\tau_j\tau_k^{-1}q^{\tilde{\nu}_j-\tilde{\nu}_k-1})}
{(1-t\tau_j^{-1}\tau_kq^{\tilde{\nu}_k-\tilde{\nu}_j})
(1-t\tau_j\tau_k^{-1}q^{\tilde{\nu}_j-\tilde{\nu}_k-1})}
\end{equation*}
(where we have used the convention that an empty product is equal to one),
which cancel each other in the final relation \eqref{relVweight}.
The verification of the difference equations
\eqref{CV+}, \eqref{CV-} is
straightforward using the explicit expressions
for $C^{\text{\tiny qR}}_+$,
$C^{\text{\tiny qR}}_-$
and $V_{\pm j}$ (in \eqref{Cqr+}-\eqref{c0} and \eqref{coefV})
and the elementary shift property
$(a;q)_{l+1}= (1-aq^l)(a;q)_l$ for the $q$-shifted factorial
($l=0,1,2,3\ldots$).
\end{proof}

\begin{lemma}\label{reslem}
For $\nu \in\Lambda_N$ \eqref{alcove}
and parameters subject to the truncation condition
\eqref{trunc},
one has that
\begin{equation}
V_{\varepsilon j}(\tau q^{\nu}) = 0\;\; \text{if}\;\;
\nu +\varepsilon e_j\not\in\Lambda_N
\;\;\;\;\; (\varepsilon = +1\;\text{or}\; -1).
\end{equation}
(Here $V_{\varepsilon j}(z)$ and $\tau q^\nu$ are taken from \eqref{coefV}
and \eqref{grid}.)
\end{lemma}
\begin{proof}
There are two situations that need to be
distinguished: either $\nu +\varepsilon e_j$ is in
$\Lambda$ \eqref{cone} but outside the
the alcove $\Lambda_N$ \eqref{alcove}, or $\nu +\varepsilon e_j$ does
not even lie in the cone $\Lambda$ \eqref{cone}.

The first situation occurs (only) when $\varepsilon =+1$
with $j=1$ and $\nu_1=N$.
We then have that $V_{+1}(\tau q^\nu)=0$ because
in the numerator the factor
$(1-t_bz_1)=(1-t_b\tau_1q^{\nu_1})=(1-t^{n-1} t_at_bq^N)$
becomes identical to
zero in view of the truncation condition \eqref{trunc}.

The second situation can occur both when $\varepsilon =+1$
or when $\varepsilon =-1$.
For $\varepsilon =+1$ we have that
$\nu+\varepsilon e_j\not\in\Lambda$
iff $j>1$ and $\nu_{j-1}=\nu_j$.
Then $V_{+j}(\tau q^\nu)=0$ because in the numerator
the factor
$(1-tz_jz_{j-1}^{-1})=(1-t\tau_j\tau^{-1}_{j-1}q^{\nu_j-\nu_{j-1}})=
(1-tt^{-1})$ is identically zero.
For $\varepsilon =-1$ we have that
$\nu+\varepsilon e_j\not\in\Lambda$ iff
either $j=n$ and $\nu_n=0$ or if
$j<n$ and $\nu_{j+1}=\nu_j$. In the former case
$V_{-n}(\tau q^\nu)=0$ because the numerator
contains a factor
$(1-t_a z_n^{-1})=(1-t_a \tau_n^{-1}q^{-\nu_n})=(1-t_at_a^{-1})=0$,
whereas in the latter case one has
that $V_{-j}(\tau q^\nu)=0$ because in the numerator one has
a factor
$(1-t z_j^{-1} z_{j+1})=(1-t \tau_j^{-1}\tau_{j+1} q^{-\nu_j+\nu_{j+1}})=
(1-t t^{-1}) =0.$
\end{proof}

After these preliminaries we are finally set
to prove the symmetry relation
\begin{equation}\label{symrel}
\langle D f , g \rangle_N^{\text{\tiny qR}} = \langle f
,D g\rangle_N^{\text{\tiny qR}} \;\;\;\;\;\;\;\;\;\;\;\;\;\;\;\;\;\;\;
(f,g\in \mathcal{H}^{\text{\tiny qR}}_N)
\end{equation}
for parameters subject to the truncation condition
\eqref{trunc}.
Evidently this amounts to showing that (cf. the definition of
$\langle \cdot ,\cdot \rangle_N^{\text{\tiny qR}}$ in \eqref{discreteIP})
\begin{eqnarray}\label{id1}
&&\sum_{\stackrel{\nu\in\Lambda_N}{1\leq j\leq n,\varepsilon =\pm 1}}
V_{\varepsilon j}(\tau q^\nu)\,
f (\tau q^{\nu+\varepsilon e_j})\,
g (\tau q^{\nu})\,
\Delta^{\text{\tiny qR}}(\nu ) = \\
&&\makebox[5em]{}
\sum_{\stackrel{\tilde{\nu}\in\Lambda_N}{1\leq j\leq n,\varepsilon =\pm 1}}
V_{-\varepsilon j}(\tau q^{\tilde{\nu}})\,
f (\tau q^{\tilde{\nu}})\,
g (\tau q^{\tilde{\nu} -\varepsilon e_j})\,
\Delta^{\text{\tiny qR}}(\tilde{\nu} ) .\nonumber
\end{eqnarray}
(At both sides of this equation
the coefficients $V_{\pm j}$
are well-defined as rational
expressions in the parameters subject to the truncation condition
\eqref{trunc},
i.e., no denominator becomes identical to zero.)
Using Lemma~\ref{reslem} the sums at both sides of \eqref{id1}
can be restricted resulting in the equation
\begin{eqnarray}\label{id2}
&&\sum_{\stackrel{1\leq j\leq n}{\varepsilon =\pm 1}}
\sum_{\stackrel{\nu\in\Lambda_N\; \text{with} }
               {\nu +\varepsilon e_j\in\Lambda_N}}
V_{\varepsilon j}(\tau q^\nu)\,
f (\tau q^{\nu+\varepsilon e_j})\,
g (\tau q^{\nu})\,
\Delta^{\text{\tiny qR}}(\nu ) = \\
&&\makebox[5em]{}
\sum_{\stackrel{1\leq j\leq n}{\varepsilon =\pm 1}}
\sum_{\stackrel{\tilde{\nu}\in\Lambda_N\; \text{with}}
               {\tilde{\nu} -\varepsilon e_j\in\Lambda_N}}
V_{-\varepsilon j}(\tau q^{\tilde{\nu}})\,
f (\tau q^{\tilde{\nu}})\,
g (\tau q^{\tilde{\nu} -\varepsilon e_j})\,
\Delta^{\text{\tiny qR}}(\tilde{\nu} ) . \nonumber
\end{eqnarray}
To check the identity \eqref{id2} (and thus proving \eqref{symrel}),
one uses Lemma~\ref{fliplem} to infer
that for given $j$ and $\varepsilon$ each term on the l.h.s. coincides
with a term on the r.h.s., where $\nu$ and $\tilde{\nu}$ are related by
$\tilde{\nu}=\nu +\varepsilon e_j$.
Phrased in other words: substituting in (the terms
corresponding to given $j$ and $\varepsilon$ at) the r.h.s.
of \eqref{id2} $\tilde{\nu}=\nu +\varepsilon e_j$ and
invoking of Lemma~\ref{fliplem} results in
(the corresponding terms at) the l.h.s. of \eqref{id2},
which completes the proof of Proposition~\ref{symprp}.

\begin{remark}
Notice that in the restriction of the sums, i.e. in passing
from \eqref{id1} to \eqref{id2}, the truncation condition
\eqref{trunc} became essential.
Moreover, it was only after this restriction of the sums that
the terms at both sides of the equation could be put into
one-to-one correspondence.
For generic parameters not satisfying any truncation condition
the proof breaks down at this step because there will be
nonzero terms of the form
\begin{equation*}
V_{+1}(\tau q^\nu)\, f (\tau q^{\nu +e_1} )\, g (\tau q^\nu )\,
\Delta^{\text{\tiny qR}}(\nu )
\end{equation*}
with $\nu_1 =N$ in the l.h.s. of \eqref{id1} that do not match with
terms in the r.h.s.
and, reversely,
will there be nonzero terms of the form
\begin{equation*}
V_{+1}(\tau q^{\tilde{\nu}})\,
f (\tau q^{\tilde{\nu}})\,
g (\tau q^{\tilde{\nu} + e_1})\,
\Delta^{\text{\tiny qR}}(\tilde{\nu} )
\end{equation*}
with $\tilde{\nu}_1=N$ that have no counter part in the l.h.s.
Therefore, in general
equation \eqref{id1} (and hence the symmetry
of $D$ with respect to $\langle \cdot ,\cdot \rangle_N^{\text{\tiny qR}}$)
will no longer hold
if no truncation condition is assumed.
\end{remark}

\section{Computation of orthonormalization constants\\
         using Pieri type formulas}\label{appb}
The key to the computation of the sum
$\langle p_\lambda ,p_\lambda \rangle^{\text{\tiny qR}}_N$
is a system of Pieri type recurrence relations
for the renormalized multivariable Askey-Wilson polynomials
\begin{equation}\label{reno}
P_\lambda (z) \equiv \hat{C}_+^{\text{qR}}(\lambda )p_\lambda(z),
\;\;\;\;\;\;\;\;\;\;\;\;\;\;\;\;\;\; \lambda \in \Lambda
\end{equation}
with $\hat{C}_+^{\text{qR}}(\lambda )$ given by \eqref{hatCqr+}.
Basically, the recurrence relations in question
provide explicit expansion
formulas of the type $E_r(z) P_\lambda (z) =\sum_\mu c_\mu P_\mu (z)$
for the products of the basis elements $P_\lambda (z)$ with
certain $W$-invariant polynomials
$E_1(z),\ldots ,E_n(z)$ that form a set
of generators for the algebra $\mathcal{H}^W$
of all $W$-invariant Laurent polynomials in the variables
$z_1,\ldots ,z_n$.
Specifically, we have (see \cite{die:self-dual,die:properties})
\begin{equation}\label{pfor}
E_r (z;\tau )\: P_{\lambda} (z)=
\sum_{\stackrel{J\subset \{ 1,\ldots ,n\} ,\, 0\leq|J|\leq r}
               {\varepsilon_j=\pm 1,\; j\in J;\;
                e_{\varepsilon J} +\lambda \in \Lambda}}
\!\!\!\!\!\!\!\!\!
\hat{U}_{J^c,\, r-|J|}(\hat{\tau}q^\lambda )\,
\hat{V}_{\varepsilon J,\, J^c}(\hat{\tau}q^\lambda )\,
P_{\lambda +e_{\varepsilon J}} (z),
\end{equation}
$r=1,\ldots ,n$, where
\begin{eqnarray*}
E_r(z;\tau )&= & \!\!\!\!
\sum_{\stackrel{ J\subset \{ 1,\ldots ,n\} }{0\leq |J|\leq r}}
\Bigl( (-1)^{r-|J|} \prod_{j\in J} (z_j+z_j^{-1}) \\
& & \makebox[4em]{}\times \!\!\!\!\!
\sum_{r\leq l_1\leq \cdots \leq l_{r-|J|}\leq n}
(\tau_{l_1}+\tau_{l_1}^{-1})\cdots (\tau_{l_{r-|J|}}+ \tau_{l_{r-|J|}}^{-1})
\Bigr)
\end{eqnarray*}
(with $\tau$, $\hat{\tau}$ taken from \eqref{tau}, \eqref{hattau})
and the expansion coefficients are governed by
\begin{eqnarray*}
\hat{V}_{\varepsilon J,\, K}(z)\!\!\! &=&\!\!\!
\prod_{j\in J} \hat{w}(z_j^{\varepsilon_j})
\prod_{\stackrel{j,k \in J}{j<k }}
\hat{v}(z_j^{\varepsilon_j}z_{k}^{\varepsilon_{k}})\,
\hat{v}(qz_j^{\varepsilon_j}z_{k}^{\varepsilon_{k}})
\prod_{\stackrel{j\in J}{k\in K}} \hat{v}(z_j^{\varepsilon_j} z_k)\,
\hat{v}(z_j^{\varepsilon_j} z_k^{-1} ),\\
\hat{U}_{K,p}(z)\!\!\! &=& \!\!\!
 (-1)^p\!\!\!\!
\sum_{\stackrel{L\subset K,\, |L|=p}
               {\varepsilon_l =\pm 1,\; l\in L }}\!\!
\Bigl( \prod_{l\in L} \hat{w}(z_l^{\varepsilon_l})
\prod_{\stackrel{l,k \in L}{l<k}}
\hat{v}(z_l^{\varepsilon_l}z_{k}^{\varepsilon_{k}})\,
\hat{v}(q^{-1}z_l^{-\varepsilon_l}z_{k}^{-\varepsilon_{k}} )\\
& &\makebox[7em]{}\times
\prod_{\stackrel{l\in L}{k\in K\setminus L}}
\hat{v}(z_l^{\varepsilon_l} z_k)\,
\hat{v}(z_l^{\varepsilon_l} z_k^{-1}) \Bigr)  ,
\end{eqnarray*}
with
\begin{equation*}
\hat{v}(\zeta )= t^{-1/2} \left( \frac{ 1-t\, \zeta }
                 { 1-\zeta} \right) ,\;\;\;\;\;\;
\hat{w}(\zeta ) =
(\hat{t}_0\hat{t}_1\hat{t}_2\hat{t}_3 q^{-1})^{-1/2}
\frac{ \prod_{0\leq r\leq 3}
(1-\hat{t}_r\zeta)}
       {(1-\zeta^2)
        (1-q\zeta^2)  } .
\end{equation*}
Here we have employed the notation
\begin{equation*}
e_{\varepsilon J} \equiv \sum_{j\in J} \varepsilon_j e_j
\end{equation*}
and we also used the conventions that empty products are equal to one,
that $\hat{U}_{K,p}\equiv 1$ for $p=0$, and
that the second sum in $E_r(z;\tau)$ is equal to one when $|J|=r$.

Even though
formula \eqref{pfor} holds for generic parameters, it is
not an entirely trivial matter to perform the reduction to the
case of parameters satisfying the truncation condition \eqref{trunc}.
The problem is
that for parameters subject to the truncation condition
the $c$-function $\hat{C}_+^{\text{\tiny qR}}(\lambda )$
\eqref{Cqr+} becomes infinite for $\lambda \in \Lambda\setminus \Lambda_N$.
(Recall that we have a zero in the denominator from the factor
$(\hat{t}_b\hat{\tau}_1;q)_{\lambda_1}=
(t^{n-1}\hat{t}_a\hat{t}_b;q)_{\lambda_1}
=(t^{n-1}t_at_b;q)_{\lambda_1}$, which is zero
for $\lambda_1>N$ when the
parameters satisfy the truncation relation $t^{n-1}t_at_b=q^{-N}$.)
Consequently, the renormalized polynomials $P_\lambda (z)$
\eqref{reno} are no longer well-defined for such parameters
when $\lambda$ lies outside the alcove $\Lambda_N$ \eqref{alcove}.
Given the fact that even if we start in the l.h.s.
with a polynomial associated
to a weight $\lambda\in\Lambda_N$, we may
end up in the r.h.s. with some polynomials
corresponding to weights outside the alcove $\Lambda_N$, it is clear
that in its present form the Pieri type recurrence formulas \eqref{pfor}
do not make sense for all $\lambda\in\Lambda_N$
when the parameters satisfy the truncation condition \eqref{trunc}.

\begin{lemma}\label{hatreslem}
For $\lambda \in\Lambda_N$ and
$\lambda +e_{\varepsilon J}\in\Lambda\setminus\Lambda_N$,
we have that
$\hat{V}_{\varepsilon J,J^c}(\hat{\tau} q^{\lambda})$
is the product of a factor of the form $(1-\hat{t}_a\hat{t}_bt^{n-1}q^{N})$
and an expression that is rational in the parameters subject
to the truncation condition \eqref{trunc} (no denominator
becomes zero after imposing \eqref{trunc}).
\end{lemma}
\begin{proof}
For $\lambda\in\Lambda_N$ one has that
$\lambda + e_{\varepsilon J}\in \Lambda\setminus \Lambda_N$
iff $\lambda_1=N$ and the index set $J\subset \{ 1,\ldots ,n\}$
contains the number $1$ with
$\varepsilon_1=+1$. Then $\hat{V}_{\varepsilon J,J^c}(\hat{\tau} q^{\lambda})$
picks up a factor $(1-\hat{t}_a\hat{t}_bt^{n-1}q^{N})$ from
the part $\hat{w}(\hat{\tau}_1 q^{\lambda_1})=\hat{w}(\hat{t}_at^{n-1}q^N)$.
\end{proof}
Lemma~\ref{hatreslem} tells us that the factor
$(1-\hat{t}_a\hat{t}_bt^{n-1}q^{N})^{-1}$
in $P_{\lambda + e_{\varepsilon J}}(z)$
when $\lambda +e_{\varepsilon J}\in \Lambda\setminus \Lambda_N$
(stemming from the denominator of the normalization factor
$\hat{C}_+^{\text{\tiny qR}}(\lambda +\varepsilon e_J )$),
is compensated in formula \eqref{pfor} by a corresponding factor
$(1-\hat{t}_a\hat{t}_bt^{n-1}q^{N})$ in the
numerator of $\hat{V}_{\varepsilon J,\, J^c}(\hat{\tau}q^\lambda)$.
Moreover, by combining this observation with Proposition~\ref{findimredprp}
it is seen that if
we restrict the variable $z$
to the grid points $\tau q^\nu$, $\nu\in\Lambda_N$, then
we end up with a recurrence formula
of the form in \eqref{pfor} in which
the sum in the r.h.s. gets restricted
to the weights of the form $\lambda + e_{\varepsilon J}$ that lie
inside the alcove $\Lambda_N$.
\begin{proposition}\label{pieriprp}
For parameters subject to the truncation condition \eqref{trunc}
and $\lambda,\nu\in\Lambda_N$ one has that
\begin{equation*}
E_r (\tau q^\nu ;\tau )\: P_{\lambda} (\tau q^\nu )=
\!\!\!\!\!\!
\sum_{\stackrel{J\subset \{ 1,\ldots ,n\} ,\, 0\leq|J|\leq r}
               {\varepsilon_j=\pm 1,\; j\in J;\;
                e_{\varepsilon J} +\lambda \in \Lambda_N}}
\!\!\!\!\!\!\!\!\!
\hat{U}_{J^c,\, r-|J|}(\hat{\tau}q^\lambda )\,
\hat{V}_{\varepsilon J,\, J^c}(\hat{\tau}q^\lambda )\,
P_{\lambda +e_{\varepsilon J}} (\tau q^\nu ),
\end{equation*}
$r=1,\ldots ,n$.
\end{proposition}
After shrinking of the domain of the variables to the grid points
in the Pieri formulas and
implementation of the truncation condition, we are now ready to compute the
the sums $\langle p_\lambda ,p_\lambda \rangle_N^{\text{\tiny qR}}$,
$\lambda \in \Lambda_N$ by means of the same method that also
led to the norms of the polynomials in the continuous case
\cite{die:self-dual}.
Replacing the products at both sides of the identity
\begin{equation}
\langle E_r P_\lambda , P_{\lambda +\omega_r}\rangle_N^{\text{\tiny qR}}
=
\langle  P_\lambda , E_r P_{\lambda +\omega_r}\rangle_N^{\text{\tiny qR}} ,
\;\;\;\;\;\;\;\;\;\;
\omega_r\equiv e_1+\cdots +e_r
\end{equation}
by the corresponding r.h.s. of the Pieri formula in
Proposition~\ref{pieriprp} and using the orthogonality of
the polynomials $P_\mu$ with respect to the bracket
$\langle \cdot ,\cdot \rangle_N^{\text{\tiny qR}}$
leads us to a relation between
$\langle P_{\lambda},P_{\lambda}\rangle_N^{\text{\tiny qR}}$ and
$\langle P_{\lambda +\omega_r},
P_{\lambda +\omega_r}\rangle_N^{\text{\tiny qR}}$ in
terms of the coefficients of the Pieri formula
\begin{eqnarray}\label{raisefund}
& &\hat{V}_{\{1,\ldots ,r\} ,\{r+1,\ldots ,n\} }(\hat{\tau} q^\lambda )
\langle P_{\lambda +\omega_r} ,P_{\lambda +\omega_r} \rangle_N^{\text{\tiny
qR}}
= \\
& & \makebox[12em]{} \hat{V}_{\{-1,\ldots ,-r\} ,\{r+1,\ldots ,n\} }
(\hat{\tau} q^{\lambda +\omega_r})
\langle P_\lambda ,P_\lambda \rangle_N^{\text{\tiny qR}}  \nonumber
\end{eqnarray}
(recall $\hat{U}_{K,p}=1$ for $p=0$).
In the above manipulations we have assumed
that $\lambda$ and $\lambda +\omega_r$ (i.e. $\lambda$ augmented
by the fundamental weight vector $\omega_r=e_1+\cdots + e_r$)
lie in the cone $\Lambda_N$ and furthermore
that the parameters satisfy the truncation condition.
To solve the recurrence relation \eqref{raisefund} for
$\langle P_{\lambda},P_{\lambda}\rangle_N^{\text{\tiny qR}}$
we exploit the following connection between the $c$-functions
$\hat{C}^{\text{\tiny qR}}_\pm$ and the coefficients
$\hat{V}_{ \{ \pm 1,\ldots ,\pm r\} ,\, \{ r+1,\ldots ,n\}}$
\begin{eqnarray*}
\frac{\hat{C}^{\text{\tiny qR}}_+(\tilde{\lambda} )}
     {C^{\text{\tiny qR}}_+(\tilde{\lambda}+\omega_r)}
&=& \hat{V}_{ \{ 1,\ldots ,r\} ,\, \{ r+1,\ldots ,n\} }
(\hat{\tau} q^{\tilde{\lambda}})\;\;\;\;\;\;\;\;\;\;\;\;\;\;\;
\text{for}\;\;\;\;\tilde{\lambda},\tilde{\lambda}+\omega_r \in \Lambda_N, \\
\frac{\hat{C}^{\text{\tiny qR}}_-(\tilde{\lambda} +\omega_r)}
     {\hat{C}^{\text{\tiny qR}}_-(\tilde{\lambda})}
&=& \hat{V}_{ \{ -1,\ldots ,-r\} ,\, \{ r+1,\ldots ,n\} }
(\hat{\tau} q^{\tilde{\lambda} +\omega_r})\;\;\;\;\;
\text{for}\;\;\;\;\tilde{\lambda}, \tilde{\lambda}+\omega_r \in \Lambda_N,
\end{eqnarray*}
which is not difficult to derive with
the aid the elementary shift property
$(a;q)_{l+1}= (1-aq^l)(a;q)_l$ for the $q$-shifted factorial
($l=0,1,2,3\ldots$).
(Notice also that for $r=1$ these two formulas
amount to the (dual versions)
of the formulas \eqref{CV+}, \eqref{CV-}
in the proof of Lemma~\ref{fliplem} specialized to the case $j=1$.)
Using these two relations we can eliminate
the coefficient functions
$\hat{V}_{ \{ \pm 1,\ldots ,\pm r\} ,\, \{ r+1,\ldots ,n\}}$
from \eqref{raisefund} entailing
\begin{equation}\label{normrel}
\langle P_\lambda ,P_\lambda \rangle_N^{\text{\tiny qR}}
\hat{\Delta}^{\text{\tiny qR}}(\lambda )  =
\langle P_{\lambda+\omega_r} ,P_{\lambda +\omega_r} \rangle_N^{\text{\tiny qR}}
\hat{\Delta}^{\text{\tiny qR}}(\lambda +\omega_r )
\end{equation}
(where
$\hat{\Delta}^{\text{\tiny qR}}=
1/(\hat{C}^{\text{\tiny qR}}_+\hat{C}^{\text{\tiny qR}}_-)$
again denotes the `Plancherel' measure.
Since the (fundamental weight) vectors $\omega_1,\ldots ,\omega_n$
positively generate the cone $\Lambda$ \eqref{cone}
it follows that the l.h.s. of \eqref{normrel}
does not depend on $\lambda\in\Lambda_N$
and so we obtain by comparing with the evaluation in $\lambda =0$
(so $P_\lambda =1$) that
\begin{equation}
\langle P_\lambda ,P_\lambda \rangle_N
= \frac{\langle 1 , 1\rangle_N}{\hat{\Delta}^{\text{\tiny qR}}(\lambda ) },
\;\;\;\;\;\;\;\;\;\; \lambda\in\Lambda_N
\end{equation}
which reads in monic form
\begin{equation}
\langle p_\lambda ,p_\lambda \rangle_N
= \mathcal{N}^{\text{\tiny qR}}(\lambda )  \langle 1 , 1\rangle_N
\;\;\;\;\;\;\;\;\;\; \lambda\in\Lambda_N
\end{equation}
with $\mathcal{N}^{\text{\tiny qR}}(\lambda )=
\hat{C}^{\text{\tiny qR}}_-(\lambda )/\hat{C}^{\text{\tiny qR}}_+(\lambda )$.

\begin{remark}
The proof for the Pieri type recurrence formulas in \cite{die:self-dual}
is complete only for parameters satisfying a self-duality
condition of the type
\begin{equation}\label{self-dual}
t_aq=t_bt_ct_d
\end{equation}
(where $t_c$ and $t_d$ denote
the two parameters complementing $t_a$ and $t_b$
such that $\{ t_a,t_b,t_c,t_d\} =\{ t_0,t_1,t_2,t_3\}$).
This condition on the parameters implies that $\hat{t}_r=t_r$
and thus that $\hat{P}_\lambda (z)=P_\lambda(z)$.
As was pointed out in Section~7.2 of
\cite{die:self-dual}, however, the Pieri type
recurrence formulas would immediately follow
for general parameters without self-duality
condition once one would
succeed in proving that Macdonald's
conjectured evaluation formula stating
that $P_\lambda (\tau )=1$ holds for such
general parameters (cf. also \cite[Theorem 3]{die:self-dual}
for a proof of the evaluation formula in the self-dual
case with parameters subject to the condition \eqref{self-dual}).
At the `CRM Workshop on algebraic methods and $q$-special functions'
in Montreal, May 1996, we learned from
Prof. Macdonald that he has
managed to produce such proof for the evaluation
formula with general parameters
using an extension of
the Cherednik approach towards the Macdonald polynomials.
(See \cite{che:double,che:macdonalds,mac:affine} for this approach,
which is deeply connected with
the representation theory of affine Hecke algebras.)
Therefore, we have formulated all our results here without imposing
the self-duality condition, even though our own direct
proof (i.e. without using the representation theory
of affine Hecke algebras) at present actually requires
at one point, viz. the verification of
the Pieri type recurrence formulas in \cite{die:self-dual}
(or equivalently the proof of the evaluation formula
$P_\lambda (\tau )=1$),
assuming that this additional condition be satisfied.
\end{remark}

\begin{remark}\label{restrem}
Another consequence of the evaluation formula
$P_\lambda (\tau )=1$ is (see \cite{die:self-dual}) the duality
relation for the
renormalized multivariable Askey-Wilson polynomials
originally conjectured by Macdonald
\begin{equation}
P_\mu (\tau q^\nu ) = \hat{P}_\nu (\hat{\tau} q^\mu ),\;\;\;\;\;\;\;\;\;
\mu ,\nu \in \Lambda ,
\end{equation}
where $\hat{P}_{\nu} (z)=C_+^{\text{\tiny qR}}(\nu)\hat{p}_\nu (z)$
denotes the renormalized multivariable Askey-Wilson polynomial
dual to $P_{\nu}(z)=\hat{C}_+^{\text{\tiny qR}}(\nu)p_\nu (z)$, i.e.
with the parameters $t_r$ being replaced by the dual parameters
$\hat{t}_r$ \eqref{hattr}.
Applying the duality relation to the restricted recurrence
formulas of Proposition~\ref{pieriprp} leads us (up to dualization)
to a system of discrete
difference equations for the multivariable
Askey-Wilson/$q$-Racah
polynomials with parameters subject to the truncation condition
\begin{equation*}
\sum_{\stackrel{J\subset \{ 1,\ldots ,n\} ,\, 0\leq|J|\leq r}
               {\varepsilon_j=\pm 1,\; j\in J;\;
                e_{\varepsilon J} +\nu \in \Lambda_N}}
\!\!\!\!\!\!\!\!\!
U_{J^c,\, r-|J|}(\tau q^\nu )\,
V_{\varepsilon J,\, J^c}(\tau q^\nu )\,
P_{\lambda } (\tau q^{\nu +e_{\varepsilon J}} ) =
E_r (\hat{\tau} q^\lambda ;\hat{\tau})\: P_{\lambda} (\tau q^\nu ),
\end{equation*}
$r=1,\ldots ,n$ (where $V_{\varepsilon J,\, J^c}$ and $U_{J^c,\, r-|J|}$
are the dual versions of $\hat{V}_{\varepsilon J,\, J^c}$ and
$\hat{U}_{J^c,\, r-|J|}$ with the parameters $t_r$
replacing $\hat{t}_r$).
This system of discrete difference equations is the restriction to the grid
$\tau q^\nu$, $\nu\in\Lambda_N$ of the system of analytic difference
equations for the multivariable Askey-Wilson polynomials
introduced in \cite{die:commuting} (see
also \cite{die:self-dual,die:properties}).
For $r=1$ the discrete
difference equation in question boils down
(after multiplication by a constant factor
$t^{n-1}(t_0\cdots t_3q^{-1})^{1/2}$)
to the restriction
to the grid points of the second order $q$-difference equation
in Section~\ref{sec2} (Eq.~\eqref{defcon2}) with
parameters subject to the truncation condition \eqref{trunc}:
\begin{equation*}
\sum_{\stackrel{1\leq j\leq n}
               {\varepsilon =\pm 1;\; \nu +\varepsilon e_j\in\Lambda_N}}
V_{\varepsilon j}(\tau q^\nu)\,
\left( p_\lambda (\tau q^{\nu+\varepsilon e_j}) -
     p_\lambda (\tau q^{\nu}) \right) =
E_{\lambda ,\lambda } p_\lambda (\tau q^{\nu}) ,
\end{equation*}
where $V_{\varepsilon j}$ and $E_{\lambda ,\lambda}$
are taken from \eqref{DAW} and \eqref{eigenvAW}.
\end{remark}

\section*{Acknowledgments}
One of us (JFvD) would like to thank Prof. K. Aomoto for
information regarding his Selberg type $q$-Jackson integrals,
Prof. A. N. Kirillov for drawing our attention to
the work of Gustafson \cite{gus:whipple's},
and Prof. I. G. Macdonald for informing us about his proof of the evaluation
formula $P_\lambda (\tau) =1$ for the multivariable Askey-Wilson polynomials
with general parameters.

\bibliographystyle{amsplain}

\end{document}